\documentclass[reqno,10pt,draft,a4paper]{amsart}
\usepackage{amsmath}
\usepackage{amssymb}
\usepackage{latexsym}
\usepackage{epsf}

\usepackage[dvips]{rotating}

\usepackage[all]{xy}
\usepackage{multicol}

\UseComputerModernTips
\CompileMatrices

\DeclareMathOperator{\Hom}{Hom}
\DeclareMathOperator{\Ext}{Ext}

\renewcommand{\ge}{\geqslant}
\renewcommand{\le}{\leqslant}

\newcommand{\C}{\mathbb{C}}

\newcommand{\N}{\mathbb{N}}

\DeclareMathOperator{\Mod}{mod}

\DeclareMathOperator{\im}{im}

\newcommand{\sS}{\mathcal{S}} 
\newcommand{\Deltag}{\Gamma}  
\newcommand{\id}{\mathrm{id}} 
\newcommand{\dd}{\delta}
\newcommand{\dl}{\delta_L}
\newcommand{\dr}{\delta_R}
\newcommand{\kl}{\kappa_L}
\newcommand{\kr}{\kappa_R}
\newcommand{\kk}{\kappa_{LR}}
\newcommand{\twoBTL}{\mathrm{2BTL}}

\newcommand{\ldescent}[1]{{\mathcal L (#1)}}
\newcommand{\rdescent}[1]{{\mathcal R (#1)}}
\newcommand\idest{i.e.,\ }
\newcommand\br{{\mathbf r}}
\newcommand\bs{{\mathbf s}}
\newcommand\bu{{\mathbf u}}
\newcommand\bv{{\mathbf v}}

\begin{document}
\theoremstyle{plain}
\numberwithin{subsection}{section}
\newtheorem{thm}{Theorem}[subsection]
\newtheorem{prop}[thm]{Proposition}
\newtheorem{cor}[thm]{Corollary}
\newtheorem{clm}[thm]{Claim}
\newtheorem{lem}[thm]{Lemma}
\newtheorem{conj}[thm]{Conjecture}
\theoremstyle{definition}
\newtheorem{defn}[thm]{Definition}
\newtheorem{rem}[thm]{Remark}
\newtheorem{eg}[thm]{Example}

\title{On the non-generic representation theory of the symplectic blob
  algebra}
\author{R. M. Green \and  P. P. Martin
\and A. E. Parker$^1$} 
\address{Department of Mathematics \\ University of Colorado \\
Campus Box 395 \\ Boulder, CO  80309-0395 \\ USA }
\email{rmg@euclid.colorado.edu} 
\address{Department of Mathematics \\ University of Leeds \\ Leeds,
  LS2 9JT \\ UK}
\email{ppmartin@maths.leeds.ac.uk}
\email{parker@maths.leeds.ac.uk}
\footnotetext[1]{Corresponding author}



\begin{abstract}
This paper   
reports   
advances in the study of the symplectic blob algebra.
We find a presentation for this algebra. 
We find a minimal poset for this as a quasi-hereditary algebra.
We discuss how to reduce
the number of parameters defining the algebra from 6 to 4 (or even 3)
without loss of representation theoretic generality.
We then find some non-semisimple specialisations by
calculating Gram determinants for certain cell modules (or standard
modules) using   
the good parametrisation defined. 
We finish by considering some quotients of 
specialisations of the symplectic blob algebra
which are isomorphic to Temperley--Lieb algebras of type $A$.
\end{abstract}

\maketitle

\section*{Introduction}\label{intro}

One of the topics considered in a previous paper \cite{gensymp} 
was a
new diagram algebra, the \emph{symplectic blob algebra}. In that paper
we 
investigated its generic representation theory and proved various
important properties of the algebra, for instance that it has a
cellular basis, that it is generically semi-simple,
and that it is a quotient of the Hecke algebra of type-$\tilde{C}$. 
The problems of determining a presentation for the algebra, and of
determining its  non-generic representation theory, were left open.
In this paper we prove (in section \ref{sect:pres})
an isomorphism with an algebra defined by a presentation, and 
begin to classify the cases when the algebra is not semi-simple. 

With the Temperley--Lieb \cite{TL} and blob algebras \cite{martsaleur}, 
the symplectic blob
algebra (or isomorphically, the affine symmetric Temperley--Lieb
algebra,
$b^{\phi}_{2n}$, also defined in \cite{gensymp}) 
belongs to an intriguing class of diagram realisations of 
Hecke algebra quotients. 
The first two
have representation theories beautifully and efficiently described in
alcove geometrical language, where the precise geometry (the
realisation of the reflection group in a weight space)
is determined
by the parameters of the algebra. 
In these first two algebras the ``good'' parametrisation appropriate to
reveal this structure is not that in which the algebras were first
described. Rather, it was discovered during efforts to put the low
rank
data on non-semisimple manifolds in parameter space in a coherent format. 
The determination of the representation theory of $b^{\phi}_{2n}$ 
in the non-semisimple cases is the next important problem in the 
programme initiated in \cite{gensymp}.

The paper is structured as follows. We first review the various
objects and notations and some of the basic properties of the
symplectic blob algebra that will be used in this paper.
This is followed by a statement and proof of a presentation for the
algebra. 
We then discuss an analogue of the blob good parametrisation,
and show how by sacrificing integrality 
we can reduce
the number of parameters defining the algebra from 6 to 4 (or even 3).
Following this, we use globalisation and localisation functors to find
a minimal labelling poset for this quasi-hereditary algebra, as a
first step to finding alcove geometry or a ``linkage principle''
(a geometrical block statement \cite{Jantzen}).
We then begin to tackle the question of non-semisimplicity by
calculating Gram determinants for certain cell modules (or standard
modules) using  
the good parametrisation above
(this complements work of De~Gier and Nichols \cite{degiernichols},
who effectively compute Gram determinants for a particular 
and distinct kind of cell module).
We finish by considering some quotients of 
specialisations of the symplectic blob algebra
which are isomorphic to Temperley--Lieb algebras of type $A$.
(Generically there is no such quotient, so these constructions provide
another way of detecting non-semisimple structure.)


\section{Review}\label{review}

We first review the objects and notations that will be used 
in this paper.

\subsection{The symplectic blob algebra}

Fix $n,m \in \N$, with $n+m$ even, and $k$ a field. 
A \emph{Brauer $(n,m)$-partition} is a partition of the set
$V \cup V'$ into pairs, where $V = \{ 1, 2, \ldots, n\}$
and  $V' = \{ 1', 2', \ldots, m'\}$.
Following Weyl \cite{weyl46}
we will think of these as \emph{Brauer $(n,m)$-diagrams}
 by taking a rectangle with
$n$ vertices labelled $1$ through to $n$ on the top and
$m$ vertices labelled $1'$ through to $m'$ on the bottom and
connecting
the two vertices $a$,  and $b$ with an arbitrary line embedded in the
plane of the rectangle, if the set $\{a,b\}$
occurs in the partition of $V \cup V'$.
For example: Take $n=m=5$ and the partition
$\{\{1,2\},\{1',3\},\{2',3'\}, \{4,5'\}, \{4',5\}\}$
and we obtain the diagram:
$$\epsfbox{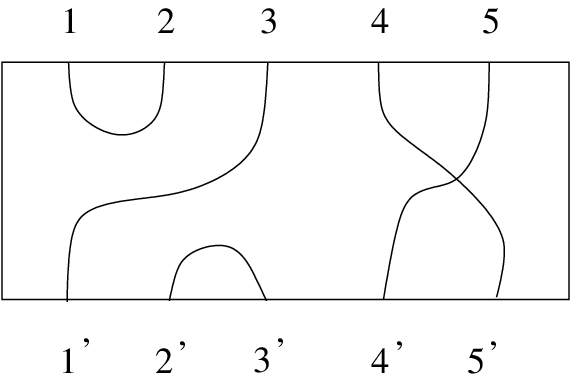}$$
Obviously
each line can be deformed isotopically 
inside the rectangle without changing the $(n,m)$-partition.
Thus any two diagrams coding 
 the same set partition are regarded as
the same diagram.

Diagrams that can be deformed isotopically within the rectangle to
obtain a diagram with no lines crossing are known as 
\emph{Temperley--Lieb diagrams}.
Thus 
$$ \epsfbox{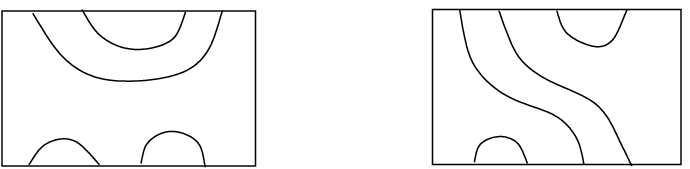}$$
are Temperley--Lieb diagrams while the 
the following diagram
$$ \epsfbox{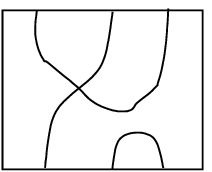}$$
is \emph{not}.
 
Our first objective is to define a certain diagram category, that is a
$k$-linear category whose hom-sets each have a basis
consisting of diagrams, and where
multiplication is defined by 
diagram concatenation, 
and 
simple straightening rules (when the concatenated object is not
formally a diagram). 
For example in our case:
$$ \epsfbox{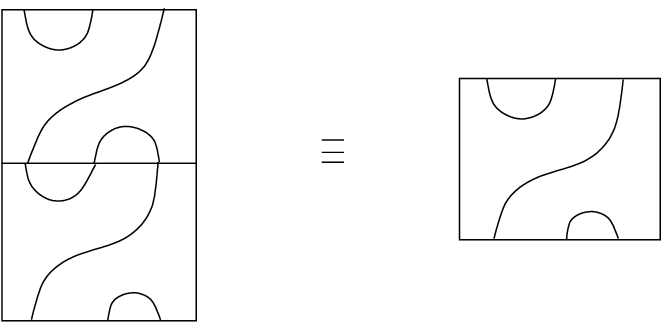}$$
by isotopy. 
When we concatenate diagrams we may get loops, for example:
\begin{equation} \label{eq:1} 
 \epsfbox{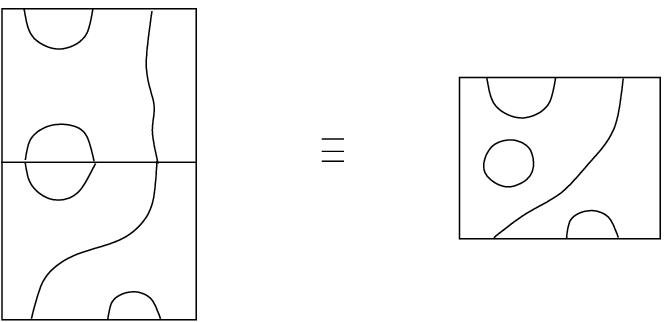}
\end{equation} 
A straightening rule is a way of expressing such products in the span
of basis diagrams. 

The resulting diagram in (\ref{eq:1}) 
is an example of a \emph{pseudo Temperley--Lieb
diagram}.
The  \emph{pseudo Temperley--Lieb diagrams} include all the
Temperley--Lieb diagrams, but we also allow diagrams with loops, which
may appear anywhere in the diagram although still with no crossing
lines. Loops that can be deformed isotopically into other loops
without crossing a line are equivalent. 
We will impose a relation on the $k$-space spanned by pseudo
Temperley--Lieb diagrams that will remove the loop.

\begin{defn}
The \emph{Temperley--Lieb algebra} $TL_n$ is the $k$-algebra
with $k$-basis the Temperley--Lieb $(n,n)$-diagrams and multiplication defined
by concatenation. We impose the relation: each loop that may arise
when multiplying is 
replaced by $\delta$, $\delta \in k$ a parameter.
\end{defn}


We can also have \emph{decorated Temperley--Lieb} diagrams where we put
elements of a monoid on the strings. When diagrams are concatenated
and words in the monoid elements form then we multiply in the monoid.
This gives us a well  defined associative 
diagram calculus --- see section 3 of \cite{gensymp}
for a detailed discussion and proof of this. Of course, we also have
\emph{decorated pseudo Temperley--Lieb} diagrams which can have
decorated loops.

We will focus on a particular set of decorated Temperley--Lieb diagrams
--- the ones that define the symplectic blob algebra.
Here our monoid that decorates the diagrams is the non-commutative
free monoid on two 
generators: a ``left'' blob, $L$, (usually a black filled-in circle 
on the diagrams) and a
``right'' blob, $R$, (usually a white filled-in circle on the diagrams).
 
A line in a (pseudo) Temperley--Lieb diagram is said to be
\emph{$L$-exposed} (respectively \emph{$R$-exposed}) if it can be deformed to touch
the left hand side (respectively right hand side) of the diagram without
crossing any other lines.

A \emph{left} (respectively \emph{right}) \emph{blob
pseudo-diagram} is a diagram obtained from a pseudo Temperley--Lieb
diagram, in which only $L$-exposed lines (respectively $R$-exposed lines)
are allowed to be decorated with  left (respectively right) blobs.
A \emph{left-right blob pseudo-diagram} is a diagram obtained
from a  pseudo Temperley--Lieb
diagram by allowing left and 
right blob decorations, with the further constraint that it must be
possible to deform decorated strings so as to take 
left blobs to the left and right blobs to the right
simultaneously.

Concatenating diagrams cannot change a $L$-exposed line to a
non-$L$-exposed line, and similarly for $R$-exposed lines. Thus the set of
left-right blob pseudo-diagrams is closed under diagram concatenation.
(See \cite[proposition 6.1.2]{gensymp}.)

The set of left-right blob pseudo-diagrams is infinite: various
features may appear. To define a finite dimensional algebra, as for
the blob algebra (see \cite[section 1.1]{marwood} for a definition)
and the Temperley--Lieb algebra (defined
above), we will 
straighten by 
replacing certain features with other features (possibly
none) multiplied by parameters from a field, $k$.

We define $B_{n,m}^{x'}$ to be the set of left-right blob
pseudo-diagrams with $n$ vertices at the top and $m$ at the bottom of the diagram
that do \emph{not} have features from the following
table.

\begin{multicols}{2}
$$
\begin{tabular}{|c|c|}
\hline
$\raisebox{-0.2cm}{\epsfbox{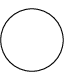}}$ 
& $\delta$\\
\hline
$\raisebox{-0.6cm}{\epsfbox{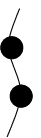}}$ 
& $\dl \raisebox{-0.6cm}{\epsfbox{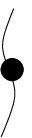}}$ \\
\hline
$\raisebox{-0.6cm}{\epsfbox{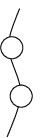}}$ 
& $\dr \raisebox{-0.6cm}{\epsfbox{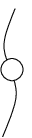}}$ \\
\hline
$\raisebox{-0.2cm}{\epsfbox{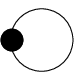}}$ 
& $\kl$\\
\hline
\end{tabular}
$$ \\
$$
\begin{tabular}{|c|c|}
\hline
$\raisebox{-0.2cm}{\epsfbox{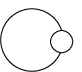}}$ 
& $\kr$\\
\hline
$\raisebox{-0.2cm}{\epsfbox{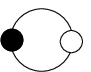}}$ 
& $\kk$\\
\hline
$\raisebox{-0.6cm}{\epsfbox{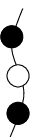}}$ 
& $k_L\raisebox{-0.6cm}{\epsfbox{lline.eps}}$ \\
\hline
$\raisebox{-0.6cm}{\epsfbox{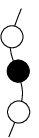}}$ 
& $k_R\raisebox{-0.6cm}{\epsfbox{rline.eps}}$ \\
\hline
\end{tabular}
$$
\end{multicols}

\begin{table}[ht]
\caption{Table encoding most of the straightening relations for $b^{x}$.\label{blobtab}}
\end{table}

The set $B_{n,m}^{x'}$ is finite and we call its elements
\emph{left-right blob diagrams}.
Define  $B_{n}^{x'}= B_{n,n}^{x'}$ and 
 $B_{}^{x'}= \bigcup_n B_{n,n}^{x'}$. 

Now define a relation on the $k$-span of all  left-right blob
pseudo-diagrams by $d \sim x d'$ if diagram $d'$ differs from $d$ 
by a substitution from left to right in the table (and extend $k$-linearly). 

A moment's thought makes it clear that to obtain a consistent set of
relations we need $RLRL= k_R RL = k_L RL$, i.e., that $k_L = k_R$.

Another (perhaps longer) moment's thought reveals that 
the $k_L$ relation is only
needed for $n$ odd and the $\kk$ relation is only needed when $n$ is
even.
It turns out to be convenient to set $\kk = k_L=k_R$.

\begin{prop}[{\cite[section 6.3]{gensymp}}]
The above relations on the $k$-span of left-right blob pseudo-diagrams 
define a finite dimensional algebra,
$b^{x'}_n$,
which has a diagram basis $B_n^{x'}$.
\end{prop}

We study this algebra by considering the quotient by the 
``topological relation'':
\begin{equation} \label{topquot}
\kappa_{LR} \;\; \raisebox{-0.2in}{\epsfbox{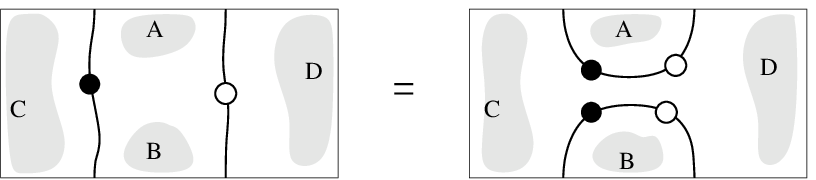}}
\end{equation}
where each labelled shaded area is shorthand for  subdiagrams that do
not have propagating lines and where
a line is called \emph{propagating} if it joins a vertex on the top
of the diagram to one on the bottom of the diagram.
(Note that there is no freedom in choosing the scalar multiple, 
once we require a  relation of this {\em form}.)

We define $B_n^x$ to be the subset of $B_n^{x'}$ that does not
contain diagrams with features as in the right hand side of relation
\eqref{topquot}.

\begin{defn}
We 
define the \emph{symplectic blob algebra}, $b_n^x$ (or $b_{n}^x(\dd,
\dl,\dr,\kl,\kr,\kk)$ if we wish to emphasise the parameters)
to be the $k$-algebra with basis $B_{n}^x$, multiplication defined
via diagram concatenation and relations as in the table above (with
$\kk=k_L=k_R$) and 
with relation \eqref{topquot}.
\end{defn}

That these relations are consistent and that we do obtain an
algebra with basis $B_n^x$ is proved in \cite[section 6.5]{gensymp}.

We have the following (implicitly assumed in \cite{gensymp}):
\begin{prop}\label{gensblob}
The symplectic blob algebra, $b_{n}^x$, is
generated by the following diagrams
\begin{multline*} 
e:=\raisebox{-0.4cm}{\epsfbox{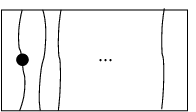}},\  
e_1:=\raisebox{-0.4cm}{\epsfbox{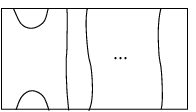}},\ 
e_2:= \raisebox{-0.4cm}{\epsfbox{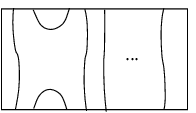}},\ 
\cdots,\ \\
e_{n-1}:=\raisebox{-0.4cm}{\epsfbox{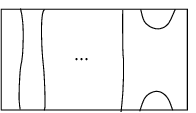}},\ 
f:=\raisebox{-0.4cm}{\epsfbox{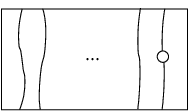}}. 
\end{multline*} 
\end{prop}
\begin{proof}
We may argue in a similar fashion as in appendix A of \cite{gensymp}
but by now inducting on the number of decorations. 
If a diagram $d$ has no decorations then the diagram is a Temperley--Lieb
diagram and the result follows.

So now assume that we have a diagram $d$ with 
$m$ decorations 
and that (for the sake of illustration) that there is a left blob --- we
would use the dual reduction in the case of a right blob.
We claim that we may use the same procedure as in the $l=0$ case
of \cite[appendix A]{gensymp}. 
If there is a decorated line starting in the first position, then we
can decompose the diagram into a product of $e$ then a diagram with one
fewer decoration. If there is no such line then take the first line
decorated with a black blob and do the same reduction as in 
\cite[appendix A]{gensymp}. 
$$\epsfbox{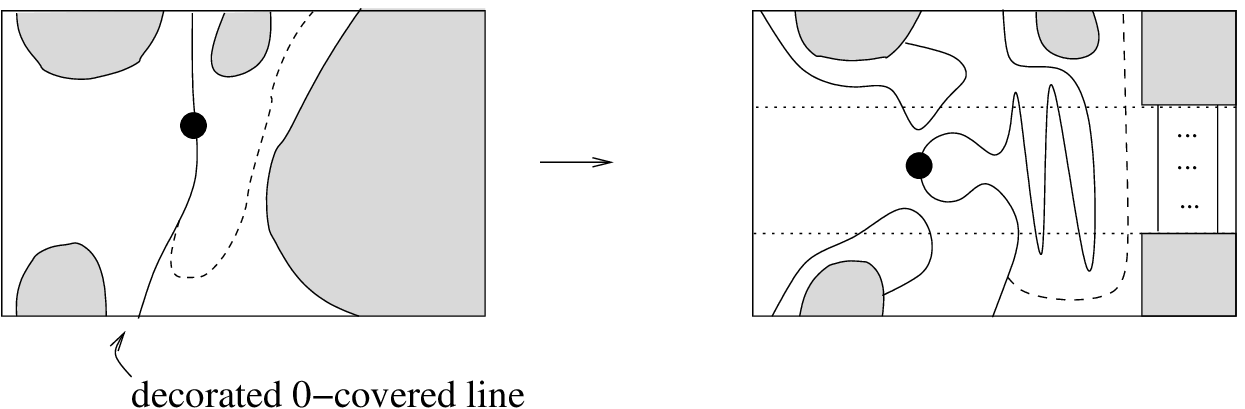}$$
$$\epsfbox{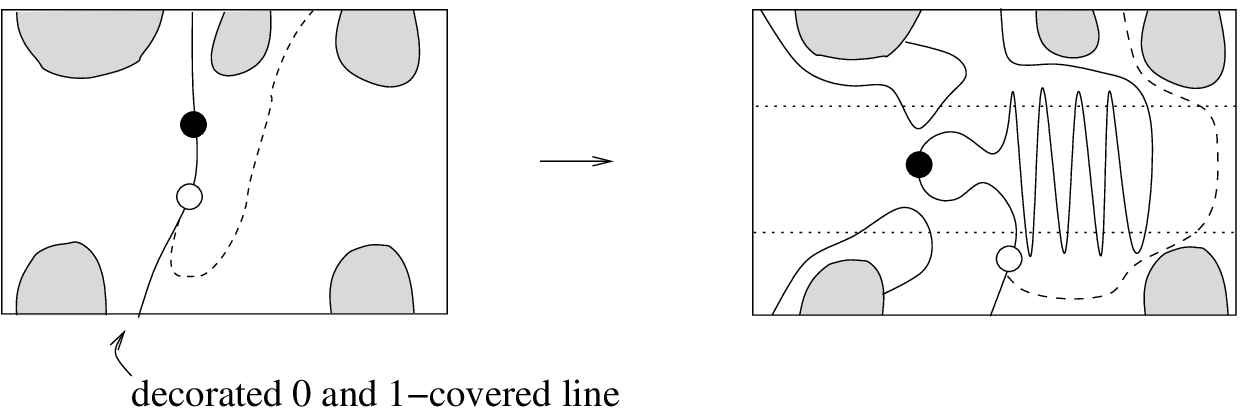}$$
The white blobs can either be moved into the shaded regions or above
or below the horizontal dotted lines. The middle region (after
``wiggling'' the line enough times)  is then the
product $e_1ee_2e_1$. The outside diagrams have strictly fewer
than $m$ decorations and hence the result follows by induction. 
\end{proof}

We now quote two results about this algebra.
\begin{prop}[{\cite[proposition 6.5.4]{gensymp}}]\label{bxeiso}
If $\dl$ is invertible then setting $e'= \frac{e}{\dl}$ we have
$$
b_{n}^x(\dd,\dl,\dr,\kl,\kr,\kk)
\cong
e'b_{n+1}^x(\dd,\kl,\dr,\dl,\kr,\kk)e'.
$$ 
\end{prop}

Similarly we have:
\begin{prop}\label{bxfiso}
If $\dr$ is invertible then setting $f'= \frac{f}{\dl}$ we have
$$
b_{n}^x(\dd,\dl,\dr,\kl,\kr,\kk)
\cong
f'b_{n+1}^x(\dd,\dl,\kr,\kl,\dr,\kk)f'.
$$ 
\end{prop}
Note the swapping of parameters required to  obtain the above
isomorphisms.

\subsection{The affine symmetric Temperley--Lieb algebra}

We now turn to defining the affine symmetric Temperley--Lieb algebra.
We may obtain (undecorated) \emph{annular Temperley--Lieb} diagrams 
by enforcing non-crossing in a cylinder or
an annulus rather than a planar rectangle:
$$ \epsfbox{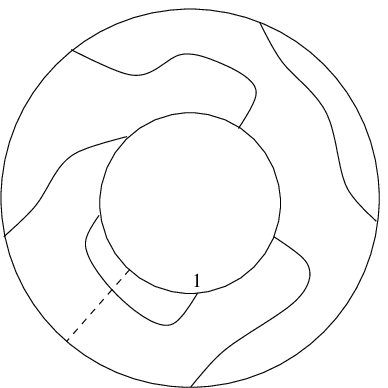}$$

We  define the annular Temperley--Lieb category by using
diagram concatenation, isotopy,
multiplication and relations of form ``loop=parameter'' as before. 
We distinguish
(contractible) loops and non-contractible
loops:
$$ \epsfbox{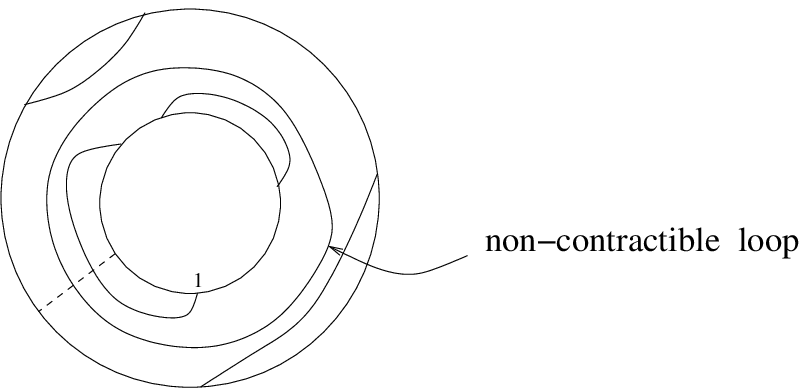}$$
Thus we will also need a relation for these non-contractible loops 
to obtain finite dimensional algebras.

(Another way to think of this is that the one loop parameter category
admits a deformation in which the loop parameter is deformed for
non-contractible loops.)

A \emph{symmetric annular Temperley--Lieb} diagram is an annular
Temperley--Lieb diagram that is symmetric about a line of symmetry ---
where we now only consider diagrams with $2n$ vertices on the top and
$2m$ on the bottom and there are $n$ (respectively $m$) vertices on 
either side of the line, for
example:
$$ \epsfbox{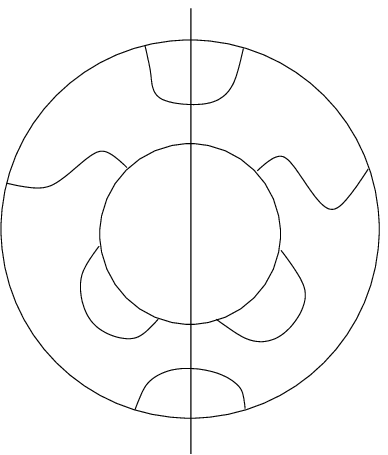}$$
The \emph{pseudo annular Temperley--Lieb} diagrams include all the
annular Temperley--Lieb diagrams, but we also allow closed loops,
including non-contractible loops.
The \emph{symmetric pseudo annular Temperley--Lieb} diagrams are the
subset of the pseudo annular Temperley--Lieb diagrams that are
symmetric.
However, we will insist on the isotopies being ``symmetric'', so that
the following two diagrams are \emph{not} equivalent:
$$ \epsfbox{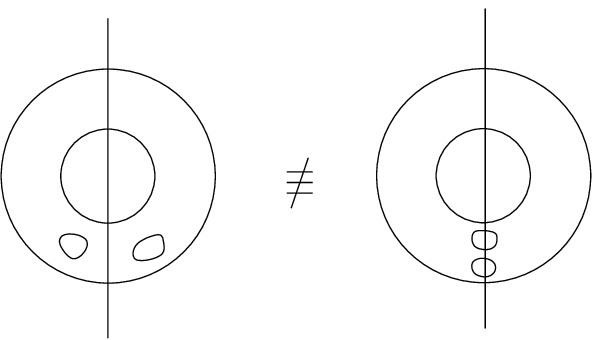}$$
This will allow us to admit a further deformation 
(see \cite[appendix B]{gensymp} for a more detailed discussion of this).

We denote the set of symmetric annular pseudo 
Temperley--Lieb diagrams by $D_{2n,2m}^\phi$, and with $n=m$  by $D_{2n}^\phi$.
(Although this is not identical with the previous paper
\cite{gensymp} ---
where $D_{2n}^\phi$ were the left-right symmetric periodic pseudo
Temperley--Lieb 
diagrams of period $2n$ --- this set is equivalent to ours
by a trivial unfolding map.)

We can (and often will) colour the annular diagrams with two colours,
black and white, such that two adjacent regions are different colours
and the marked corner (top of the 0-reflection
line) is always white. 
We also split the line of symmetry into two parts --- the 1-reflection
line and 0-reflection line as marked:
$$ \epsfbox{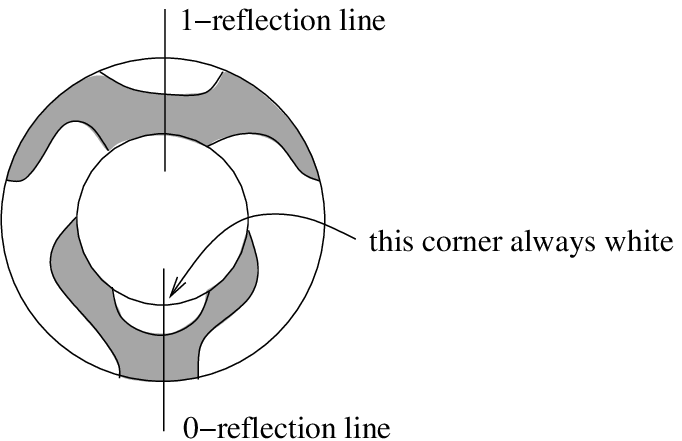}$$
If both the top and the bottom corner of the 0-reflection line are
white then we say the pseudo diagrams are \emph{colouring composable} (or
CC for short). (The above diagram is not colouring composable.)

The subset of $D_{2n}^\phi$ consisting of CC diagrams will be denoted
by $CC_{2n}$.

The set $CC_{2n}$ is closed under diagram concatenation as the
following example illustrates:
\begin{eg}
$$ \epsfbox{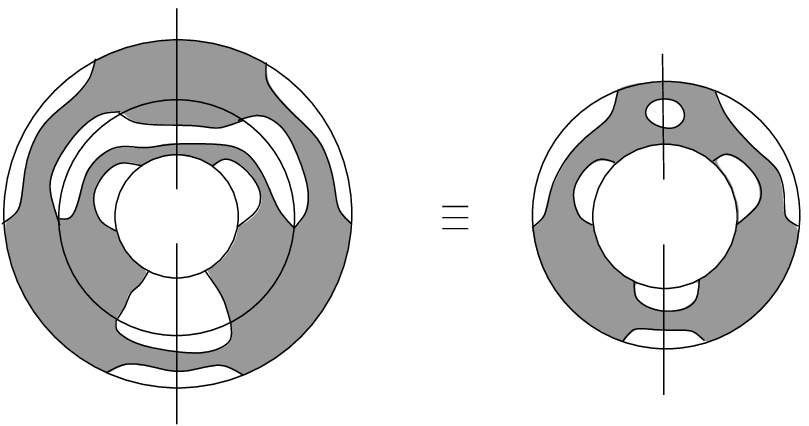}$$
\end{eg}

The set $CC_{2n}$ is not finite, so to produce a finite dimensional
diagram algebra we will need some relations.

We take $B_{2n}^\phi \subset CC_{2n}$ to be those elements of
$D_{2n}^\phi$ that do \emph{not} have features on the LHS of
table~\ref{blobtab2}. The relations will be that the features on the LHS
are replaced by the parameter on the RHS.
\begin{table}
$$
\begin{tabular}{|c|c|}
\hline
symmetric pair of contractible loops & $\delta$\\
\hline
black loop astride 0-reflection line & $\kappa_L$\\
\hline
white loop astride 0-reflection line & $\delta_L$\\
\hline
black loop astride 1-reflection line 
& $\begin{cases} \kappa_R &\mbox{if $n$  is even}\\
 \delta_R &\mbox{if $n$  is odd} \end{cases}$ \\
\hline
white loop astride 1-reflection line 
& $\begin{cases} \delta_R &\mbox{if $n$  is even}\\
 \kappa_R&\mbox{if $n$  is odd} \end{cases}$ \\
\hline
pair of non-contractible loops & $\kappa_{LR}$\\
\hline
\end{tabular}
$$
\caption{Table encoding most of the straightening relations for $b^{\phi}$.\label{blobtab2}}
\end{table}

The set $B_{2n}^\phi$ is finite. 

\begin{defn}
The \emph{affine
symmetric Temperley--Lieb algebra}, $b_{2n}^\phi$
is the $k$-algebra with basis $B_{2n}^\phi$, multiplication defined
via diagram concatenation and relations as in the table above.
\end{defn}

The reason why we introduce an odd-even dependency of the parameters for
the affine symmetric Temperley--Lieb algebra is the following
result:

\begin{prop}[{\cite[proposition 7.2.4]{gensymp}}]
The symplectic blob algebra and the affine symmetric Temperley--Lieb
algebra are isomorphic, with the obvious identification of parameters.
$$b_n^x(\dd, \dl,\dr,\kl,\kr,\kk)
\cong
b_{2n}^{\phi}(\dd, \dl,\dr,\kl,\kr,\kk)
.$$
\end{prop}

We do not stress the explicit isomorphism here as it is not needed in
the sequel. But 
as a corollary we obtain the following:
\begin{prop}\label{gens}
The affine symmetric Temperley--Lieb algebra, $b_{2n}^\phi$, is
generated by the following diagrams:
\begin{multline*} 
e:=\raisebox{-1.0cm}{\epsfbox{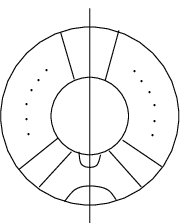}},\  
e_1:=\raisebox{-1.0cm}{\epsfbox{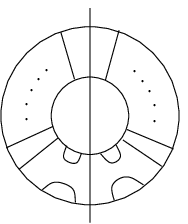}},\ 
e_2:= \raisebox{-1.0cm}{\epsfbox{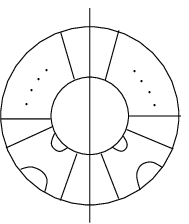}},\ 
\cdots,\ \\
e_{n-1}:=\raisebox{-1.0cm}{\epsfbox{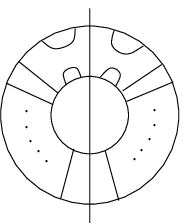}},\ 
f:=\raisebox{-1.0cm}{\epsfbox{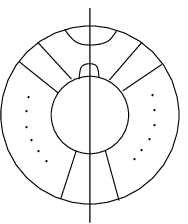}}. 
\end{multline*} 
\end{prop}
Here we have abused notation slightly and called the generators by
the same names on both sides of the isomorphism. This above identification
would also
define the isomorphism from $b_n^x$ to $b_{2n}^\phi$ if we knew a
priori that $b_{2n}^\phi$ was generated by the images of the
generators for $b_n^x$.




\section{Presenting the symplectic blob algebra}\label{sect:pres}
\newcommand{\TBTL}{two boundary Temperley--Lieb algebra}
\newcommand{\sblob}{symplectic blob algebra}

The \TBTL\ (2BTL) is (a certain parametrisation of) an
infinite-dimensional quotient of 
the Temperley--Lieb algebra of type-$\tilde C$. 
In this section we show that the \sblob\ has a presentation 
consisting of the relations for the 
 2BTL, 
together with two additional relations that make it finite
dimensional.

Presently little is known, even in the generic case, about the
larger algebras, so attacking non-trivial but tractable finite dimensional
quotients is an effective approach to their study.
Our result shows that the \sblob\ is a good tool for this study. 

We start by giving the definition of 2BTL in a suitable form. 
We must then assemble a large number of preparatory lemmas,
before finally approaching the proof of the theorem. 

\subsection{A presentation}

\begin{defn}\label{seca5}
Fix $n$. 
Let $S = \{E_0, E_1, \ldots, E_n\}$, and let $S^*$ be the free monoid on $S$.
Define the \emph{commutation monoid} $M$ to be the quotient of $S^*$ by the
relations $$
E_i E_j \equiv E_j E_i \text{ for all } 0 \leq i, j \leq n \text{ with } 
|i - j| > 1
.$$
\end{defn}

\begin{defn}\label{2BTL pres}
Let $A_n$ be the quotient of the $k$-monoid-algebra of $M$  
by the following relations:
\begin{align*}
E_0^2 &= \delta_L E_0 
&
E_1E_0E_1 &= \kappa_L E_1
\\
E_i^2 &= \delta E_i\quad\mbox{for }1 \le i \le n-1
&
E_iE_{i+1}E_i &= E_i \quad\mbox{for }1 \le i \le n-2 
\\
E_{n}^2 &= \delta_R E_n 
&
E_{i+1}E_{i}E_{i+1} &= E_{i+1}\quad\mbox{for }1 \le i \le n-2
\\
&
&
E_{n-1}E_nE_{n-1} &= \kappa_R E_{n-1}
\\
IJI &= \kappa_{LR} I 
&
JIJ &= \kappa_{LR} J 
\end{align*}
where 
$$I= \begin{cases} 
E_1 E_3 \cdots E_{2m-1} &\mbox{if $n=2m$}\\
E_1 E_3 \cdots E_{2m-1} E_{2m+1} &\mbox{if $n=2m +1$}\\
\end{cases}$$
$$J= \begin{cases} 
E_0 E_2 \cdots E_{2m-2} E_{2m} &\mbox{if $n=2m$}\\
E_0 E_2 \cdots E_{2m}  &\mbox{if $n=2m +1$}\\
\end{cases}$$
\end{defn}
Note $I=E_1$ and $J=E_0$ if $n=1$.
We will sometimes write  $E$ for $E_0$ and $F$ for $E_n$.

The algebra $A_n$ is the quotient of the two boundary Temperley--Lieb
algebra (or the Temperley--Lieb algebra of type $\tilde{C}$) by
 the additional relations  $IJI=\kk I$ and $JIJ=\kk J$.

\begin{thm}\label{thmmain}
The symplectic blob algebra 
$b_{n}^x$   
is isomorphic to the algebra
$A_n$ via an isomorphism 
$$
\phi: A_n \rightarrow b_{n}^x
$$ 
induced by $E\mapsto e$, $E_1\mapsto e_1$,
$\ldots$, $E_{n-1}\mapsto e_{n-1}$ and $F\mapsto f$.
\end{thm}
It is straightforward to check that the generators already given for 
both the affine symmetric Temperley--Lieb algebra and the symplectic
blob algebra satisfy the $A_n$ relations. Thus the map $\phi$ in the
theorem is a surjective homomorphism and hence we 
need only to prove injectivity.
The rest of this section is devoted to proving this theorem.


\subsection{Definitions associated to the monoid $M$}

Two monomials $\bu,\bu'$ in the generators $S$ are said to be \emph{commutation equivalent}
if $\bu\equiv\bu'$  
in $M$. 
The {\it commutation class}, $\overline{\bu}$, of a 
monomial $\bu$ consists of the monomials that are commutation equivalent to it.

The \emph{left descent set} (respectively, \emph{right descent set})
of a monomial $\bu$
consists of all the initial (respectively, terminal) letters of the 
elements of $\overline{\bu}$.  We denote these sets by $\ldescent{\bu}$ and
$\rdescent{\bu}$, respectively.


\begin{defn}\label{seca1}
A \emph{reduced monomial} is a monomial $\bu$ in the generators $S$
such that no $\bu' \in \overline{\bu}$ can be expressed
as a scalar multiple of a strictly shorter monomial
using the relations in Definition~\ref{2BTL pres}. 
\\
If we have $\bu = \bu_1 s \bu_2 s \bu_3$ for some generator
$s$, then the occurrences of $s$ in $\bu$ are said to be \emph{consecutive} if
$\bu_2$ contains no occurrence of $s$.
\end{defn}

\begin{defn}\label{seca2}
Two monomials in the generators, $\bu$ and $\bu'$, are said to be
\emph{weakly equivalent} if $\bu$ can be transformed into a nonzero multiple of
$\bu'$ by applying finitely many relations.  
\\
In this situation, we also say
that $D$ and $D'$ are weakly equivalent, where $D$ and $D'$ are the diagrams
equal to $\phi(\bu)$ and $\phi(\bu')$, respectively.  
If $P$ is a property that 
diagrams may or may not possess, then we say $P$ is \emph{invariant under weak
equivalence} if, whenever $D$ and $D'$ are weakly equivalent diagrams, then
$D$ has $P$ if and only if $D'$ has $P$.
\end{defn}

\begin{defn}\label{seca3}
Let $D$ be a diagram.  For $g \in \{L, R\}$ and $$
k \in \{1, \ldots, n, 1', \ldots, n'\}
,$$ we say that $D$ is \emph{$g$-decorated at the point $k$} if (a) the edge
$x$ connected to $k$ has a decoration of type $g$, and (b) the decoration of
$x$ mentioned in (a) is closer to point $k$ than any other decoration on $x$.
\end{defn}

In the sequel, we will sometimes invoke Lemma \ref{seca4} without explicit comment.

\begin{lem}\label{seca4}
The following properties of diagrams are invariant under weak equivalence:
\begin{enumerate}
\item[(i)]{the property of being $L$-decorated at the point $k$;}
\item[(ii)]{the property of being $R$-decorated at the point $k$;}
\item[(iii)]{for fixed $1 \leq i < n$, the property of points $i$ 
and $(i+1)$ being connected by an undecorated edge;}
\item[(iv)]{for fixed $1 \leq i < n$, the property of points $i'$ 
and $(i+1)'$ being connected by an undecorated edge.}
\end{enumerate}
\end{lem}

\begin{proof}
It is enough to check that each of these properties is respected by each type
of diagrammatic reduction, because the diagrammatic algebra is a homomorphic
image of the algebra given by the monomial presentation.
This presents no problems, but notice that
the term ``undecorated'' cannot be removed from parts (iii) and (iv), because
of the topological relation.
\end{proof}

Elements of the commutation monoid 
$M$ have the following normal form, established in \cite{CF}.

\begin{prop}[Cartier--Foata normal form]\label{seca6} 
Let $\bs$ be an element of the commutation monoid $M$.  Then $\bs$ has
a unique factorization in $M$ of the form $$
\bs = \bs_1 \bs_2 \cdots \bs_p
$$ such that each $\bs_i$ is a product of distinct commuting elements of
$S$, and such that for each $1 \leq j < p$ and each generator $t \in S$
occurring in $\bs_{j+1}$, there is a generator $s \in S$ occurring in
$\bs_j$ such that $st \ne ts$.
\end{prop}

\begin{rem}\label{seca7}
The Cartier--Foata normal form may be defined inductively, as follows.  
Let $\bs_1$ be the product of
the elements in $\ldescent{\bs}$.  Since $M$ is a cancellative monoid,
there is a unique element $\bs' \in M$ with $\bs = \bs_1 \bs'$.  If $$
\bs' = \bs_2 \cdots \bs_p
$$ is the Cartier--Foata normal form of $\bs'$, then $$
\bs_1 \bs_2 \cdots \bs_p
$$ is the Cartier--Foata normal form of $\bs$.
\end{rem}

\begin{rem}\label{seca8}
The monoid $M$ is useful for our purposes because the symplectic
blob algebra is a quotient of the monoid algebra of $M$, where we identify
$E_0$ with $e$, $E_i$ with $e_i$ and $E_n$ with $f$, as usual.
\end{rem}

\begin{defn}\label{secc1}
Let $\bu$ be a reduced monomial in the generators $E_0, \ldots, E_n$.  We say
that $\bu$ is \emph{left reducible} (respectively, \emph{right reducible})
if it is commutation equivalent to a monomial of the form $\bu' = st \bv$ 
(respectively, $\bu' = \bv ts$), where $s$ and $t$ are noncommuting
generators and $t \not\in \{E, F\}$.  In this situation, we say that
$\bu$ is left (respectively, right) reducible via $s$ to $t \bv$ (respectively,
to $\bv t$).
\end{defn}

\subsection{Preparatory lemmas}   


The following result is similar to \cite[Lemma 5.3]{G35}, but we give a complete argument here because the proof in
\cite{G35} contains a mistake (we thank D.C. Ernst for pointing this out).

\begin{lem}\label{secc2}
Suppose that $\bs \in M$ corresponds to a reduced monomial, 
and let $\bs_1 \bs_2 \cdots \bs_p$ be the Cartier--Foata
normal form of $\bs$.  Suppose also that $\bs$ is not left reducible.
Then, for $1 \leq i < p$ and $0 \leq j \leq n$, the following hold:
\begin{enumerate}
\item[(i)]
{if $E_0$ occurs in $\bs_{i+1}$, then $E_1$ occurs in $\bs_i$;}
\item[(ii)]
{if $E_n$ occurs in $\bs_{i+1}$, then $E_{n-1}$ occurs in $\bs_i$;}
\item[(iii)]
{if $j \not\in \{0, n\}$ and $E_j$ occurs in $\bs_{i+1}$, then both 
$E_{j-1}$ and $E_{j+1}$ occur in $\bs_i$.}
\end{enumerate}
\end{lem}

\begin{proof}
The assertions of (i) and (ii) are immediate from properties of the normal
form, because $E_1$ (respectively, $E_{n-1}$) is the only generator not
commuting with $E_0$ (respectively, $E_n$).  We will now prove (iii) 
by induction on $i$.  Suppose first that $i = 1$.

Suppose that $j \not\in \{0, n\}$ and that $E_j$ occurs in $\bs_2$.  
By definition of the normal form, there must be a generator $s \in \bs_1$ 
not commuting with $E_j$.   Now $s$ cannot be the only such generator, or 
$\bs$ would be left reducible via $s$.
Since the only generators not
commuting with $E_j$ are $E_{j-1}$ and $E_{j+1}$, these must both
occur in $\bs_1$.

Suppose now that the statement is known to be true for $i < N$, and
let $i = N \geq 2$.
Suppose also that $j \not\in \{0, n\}$ and that $E_j$ occurs in $\bs_{N+1}$.  
As in the base case, there must be at least one generator $s$ occurring 
in $\bs_N$ that does not commute with $E_j$.  

Let us first consider the case where $j \not\in \{1, n-1\}$, and write
$s = E_k$ for some $1 \leq k \leq n-1$.  The restrictions on $j$ mean that
we cannot have $E_j E_k E_j$ occurring as a subword of any reduced monomial.
However, $E_j$ occurs in $\bs_{N-1}$ by the inductive hypothesis, and this
is only possible if there is another generator, $s'$, in $\bs_N$ that does not 
commute with $E_j$.  This implies that $\{s', E_k\} = \{E_{j-1}, E_{j+1}\}$,
as required.

Now suppose that $j = 1$ (the case $j = n-1$ follows by a symmetrical
argument).  If both $E_0$ and $E_2$ occur in $\bs_N$, then
we are done.  If $E_2$ occurs in $\bs_N$ but $E_0$ does not, then the
argument of the previous paragraph applies.  Suppose then that $E_0$ occurs
in $\bs_N$ but $E_2$ does not.  By statement (i), $E_1$ occurs in 
$\bs_{N-1}$, but arguing as in the previous paragraph, we find this cannot
happen, because it would imply that $\bs$ was commutation equivalent to
a monomial of the form $\bv' E_1 E_0 E_1 \bv''$, which is incompatible
with $\bs$ being reduced.  This completes the inductive step.
\end{proof}

The following is a key structural property of reduced monomials.

\begin{prop}\label{secc3}
Suppose that $\bs \in M$ corresponds to a reduced monomial, 
and let $\bs_1 \bs_2 \cdots \bs_p$ be the Cartier--Foata
normal form of $\bs$, where $\bs_p$ is nonempty.
Suppose also that $\bs$ is neither left reducible nor right reducible.
Then either {\rm{(i)}} $p = 1$, meaning that $\bs$ is a product of commuting
generators or {\rm{(ii)}} $p = 2$ and either $\bs = IJ$ or $\bs = JI$.
\end{prop}

\begin{proof}
If $p = 1$, then case (i) must hold, so we will assume
that $p > 1$.  

A consequence of Lemma \ref{secc2} is that 
if $\bs_{i+1} = I$ then $\bs_i = J$, and 
if $\bs_{i+1} = J$ then $\bs_i = I$.
It follows that if $\bs_p \in \{I, J\}$ (as algebra elements), then 
$\bs$ must be an alternating
product of $I$ and $J$.  Since $\bs$ is reduced, this forces $p = 2$ and
either $\bs = IJ$ or $\bs = JI$.  We may therefore assume that $\bs_p \not\in
\{I, J\}$.

Since $\bs_p \not \in \{I, J\}$ and $\bs_p$ is a product of commuting 
generators, at least one of the following two situations must occur.
\begin{enumerate}
\item[(a)]{For some $2 \leq i \leq n$, $\bs_p$ contains an occurrence of $E_i$
but not an occurrence of $E_{i-2}$.}
\item[(b)]{For some $0 \leq i \leq n-2$, $\bs_p$ contains an occurrence of 
$E_i$ but not an occurrence of $E_{i+2}$.}
\end{enumerate}

Suppose we are in case (a).  In this case, Lemma \ref{secc2} 
means that there must
be an occurrence of $E_{i-1}$ in $\bs_{p-1}$; furthermore, $E_{i-1} \not\in
\{E, F\}$, because $e_{i-1}$ fails to commute with two other generators ($E_i$
and $E_{i-2}$).  However,
one of these generators, $E_{i-2}$ does not occur in $\bs_p$.  It follows 
that $\bs$ is right reducible (via $E_i$), which is a contradiction.
Case (b) leads to a similar contradiction, again involving right
reducibility, which completes the proof.
\end{proof}


\begin{lem}\label{secb1}
Let $\bu = \bu_1 s \bu_2 s \bu_3$ be a reduced word in which the occurrences
of the generator $s$ are consecutive, and suppose that every generator
in $\bu_2$ not commuting with $s$ is of the same type, $t$ say.  Then 
$\bu_2$ contains only one occurrence of $t$, and $s \in \{E, F\}$.
\end{lem}

\begin{proof}
The proof is by induction on the length, $l$, of the word $\bu_2$.  Note
that $\bu_2$ must contain at least one generator not commuting with $s$,
or after commutations, we could produce a subword of the form $ss$.
This means that the case $l = 0$ cannot occur.  

If $\bu_2$ contains only one generator not commuting with $s$, then after
commutations, $\bu$ contains a subword of the
form $sts$.  This is only possible if $s \in \{E, F\}$, and this establishes
the case $l = 1$ as a special case.

Suppose now that $l > 1$.  By the above paragraph, we may reduce to the case
where $\bu_2 = \bu_4 t \bu_5 t \bu_6$, and the indicated
occurrences of $t$ are consecutive.  Since $\bu_5$ is shorter than
$\bu_2$, we can apply the inductive hypothesis to show that $t \in \{E, F\}$
and $\bu_5$ contains only one generator, $u$, that does not commute with $t$.
We cannot have $u = s$, or the original occurrences of $s$ would not be
consecutive.  This means that $t$ fails to commute with two different
generators, contradicting the fact that $t \in \{E, F\}$ and completing
the proof.
\end{proof}

\begin{lem}\label{secb2}
Let $\bu$ be a reduced monomial.
\begin{enumerate}
\item[(i)]{Between any two consecutive occurrences of $E$ in $\bu$,
there is precisely one letter not commuting with $E$ (\idest an occurrence
of $E_1$).}
\item[(ii)]{Between any two consecutive occurrences of $F$ in $\bu$,
there is precisely one letter not commuting with $F$ (\idest an occurrence
of $E_{n-1}$).}
\item[(iii)]{Between any two consecutive occurrences of $E_i$ in $\bu$,
there are precisely two letters not commuting with $E_i$, and they correspond
to distinct generators.}
\end{enumerate}
\end{lem}

\begin{proof}
To prove (i), we apply Lemma \ref{secb1} with $s = E$; the hypotheses
are satisfied as we necessarily have $t = E_1$.  The proof of (ii) is similar.

To prove (iii), write $\bu = \bu_1 s \bu_2 s \bu_3$ for consecutive
occurrences of the generator $s = E_i$.  Since $s \not\in \{E, F\}$, the 
hypotheses of Lemma \ref{secb1} cannot be satisfied, so $\bu_2$ must have at
least one occurrence of each of $t_1 = E_{i-1}$ and $t_2 = E_{i+1}$.
Suppose that $\bu_2$ contains two or more occurrences of $t_1$.  The fact
that the occurrences of $s$ are consecutive means that two consecutive
occurrences of $t_1$ cannot have an occurrence of $s$ between them.
Applying Lemma \ref{secb1}, this means that there is precisely one generator
$u$ between the consecutive occurrences of $t_1$ such that $t_1 u \ne u t_1$,
and furthermore, that $t_1 \in \{E, F\}$.  This is a contradiction, because
$t_1$ fails to commute with two different generators ($s$ and $u$).

One can show similarly that $\bu_2$ cannot contain two or more 
occurrences of $t_2$.  We conclude that each of $t_1$ and $t_2$ occurs
precisely once, as required.
\end{proof}

\subsection{The map $\phi$}

\begin{figure}[ht]
\epsfbox{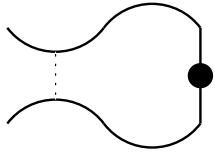}
\caption{Direction reversal of an arc}\label{fig1}
\end{figure}

\begin{lem}\label{secb3}
Let $D$ be a diagram representing a reduced monomial (\idest it is
$\phi$ applied to a reduced monomial).  Then the only way
an arc of $D$ may change direction from left to right or vice versa is
as shown in Figure \ref{fig1} or its mirror image: the turn is performed with
a vertical decorated section, and the turn is tight in the sense that the
two horizontal arcs connected by a thin dotted line in the diagram arise from 
the same letter of the monomial.
\end{lem}

\begin{proof}
This is a restatement of Lemma \ref{secb2}~(iii).
\end{proof}

\begin{lem}\label{secb4}
Let $D$ be a diagram representing a reduced monomial $\bu$ (\idest
$D=\phi(\bu)$).
\begin{enumerate}
\item[(i)]{The diagram $D$ is $L$-decorated at $1$ (respectively, $1'$)
if and only if 
the left (respectively, right) descent set of $\bu$ contains $E$.}
\item[(ii)]{The diagram $D$ is $R$-decorated at $n$ (respectively, $n'$), 
if and only if 
the left (respectively, right) descent set of $\bu$ contains $F$.}
\item[(iii)]{Suppose that $1 \leq i < n$.  Then points $i$ and $i+1$ 
(respectively, $i'$ and $(i+1)'$) in $D$ are connected by an undecorated 
edge if and only if the left
(respectively, right) descent set of $\bu$ contains $E_i$.}
\end{enumerate}
\end{lem}

\begin{proof}
In all three cases, the ``if'' statements follow easily from diagram calculus
considerations, so we only prove the ``only if'' statements.

Suppose for a contradiction that $D$ is $L$-decorated at $1$, but that the
left descent set of $\bu$ does not contain $E$.  For this to happen, the 
arc leaving point $1$ must eventually encounter an $L$-decoration, but must
first encounter a horizontal arc corresponding to an occurrence of the
generator $E_1$.  The only way this can happen and be consistent with Lemma
\ref{secb3} is for the 
arc to then travel to the eastern edge after encountering
$E_1$, then change direction and then travel back to the western edge, as shown
in Figure \ref{fig2}.  (Note that this can only happen if $n$ is odd, and that
as before, the thin dotted lines indicate pairs of horizontal edges
that correspond to the same generator.)

\begin{figure}
        \epsfbox{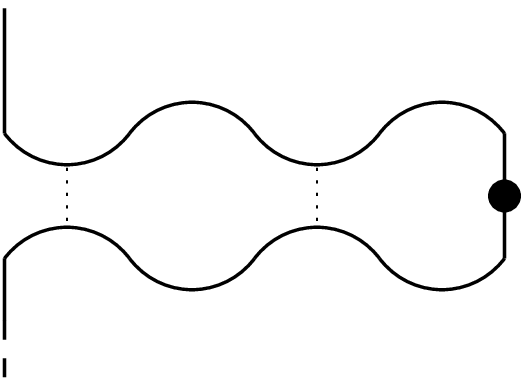}
\caption{Illustrating the proof of Lemma \ref{secb4}~(i).\label{fig2}}
\end{figure}

By Lemma \ref{seca4}, this implies that the arc connected to point $1$ is 
$R$-decorated, which in turn means that it cannot also be $L$-decorated.  This
proves that if $D$ is $L$-decorated at $1$ then the left descent set
of $\bu$ contains $E$, and the claim regarding $1'$ and the right descent
set is proved similarly.  This completes the proof of (i), and the proof of
(ii) follows by modifying the above proof in the obvious way.

We now turn to (iii).
Suppose for a contradiction that points $i$ and $i+1$ in $D$ are connected
by an undecorated edge, but that the
left descent set of $\bu$ does not contain $E_i$.  
For this to happen, it must be the case that either
(a) it is not the case that the arc leaving point $i$ 
encounters the northernmost occurrence of $E_i$ before any other generator 
or 
(b) it is not the case that the arc leaving point $i+1$ 
encounters the northernmost occurrence of $E_i$ before any other generator.
(We allow the possibility that (a) and (b) could both occur.
Notice that since the
arc crosses the line $x = i + 1/2$, it must encounter a generator $E_i$ at
some stage.)  We deal with case (a); the treatment of case (b) follows by
a symmetrical argument.  The only way case (a) can occur consistently
with Lemma \ref{secb3} is for the situation in Figure \ref{fig3}
to occur, and even this
is impossible unless $i$ is odd.

\begin{figure}
\caption{Illustrating the proof of Lemma \ref{secb4}~(iii).\label{fig3}}
        \epsfbox{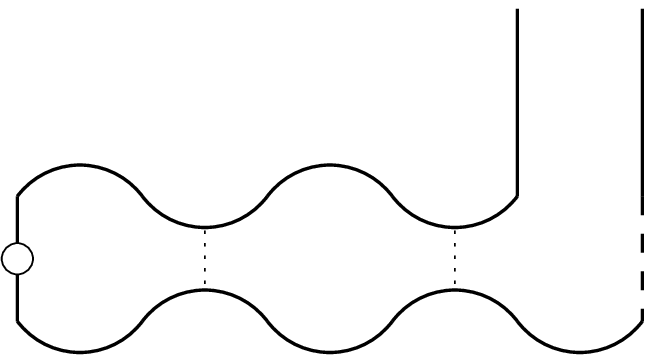}
\end{figure}

This is a contradiction by Lemma \ref{seca4}, because it implies that the arc
connected to point $i$ is $L$-decorated, and we assumed that it was 
undecorated.
This proves that if points $i$ and $i+1$ of $D$ are connected by an undecorated
edge, then the left descent set
of $\bu$ contains $E_i$, and the claim regarding the right descent
set is proved similarly.  This completes the proof of (iii).
\end{proof}

\begin{lem}\label{secb5}
Let $\bu$ and $\bu'$ be reduced monomials that map to the same
diagram $D$ under $\phi$.
\begin{enumerate}
\item[(i)]{If $\bu'$ is a product of commuting generators, then
$\bu$ and $\bu'$ are equal as algebra elements.}
\item[(ii)]{If $\bu' = IJ$ or $\bu' = JI$, then $\bu$ and $\bu'$
are equal as algebra elements.}
\end{enumerate}
\end{lem}

\begin{proof}
We first prove (i).  By Lemma \ref{secb4}, we must have $$
\ldescent{\bu} = \rdescent{\bu} = \ldescent{\bu'} = \rdescent{\bu'}
,$$ because $\bu$ and $\bu'$ are
represented by the same diagram.  

Suppose that $\bu'$ contains an occurrence of the generator $E$.
This implies (a) that $\bu$ must contain an occurrence of $E$, because
$\bu$ and $\bu'$ are represented by the same diagram, and (b) 
$E \in \ldescent{\bu} \cap \rdescent{\bu}$, by Lemma \ref{secb4}.
Suppose also (for a contradiction) that $\bu$ contains two occurrences of
the generator $E = E_0$.  By Lemma \ref{secb2}, there must be an occurrence of
$E_1$ between the first (\idest leftmost or northernmost) two occurrences
of $E_0$.  

Since points $1$ and $1'$ of $D$ are connected by an $L$-decorated edge,
there must be an occurrence of $E_2$ immediately above the aforementioned 
occurrence of $E_1$ in order to prevent the edge emerging from point $1$ 
from exiting the box at point $2$.  (``Immediately above'' means that there
are no other occurrences of $E_1$ or $E_2$ between the two occurrences
mentioned.)  In turn, we must have an occurrence of $E_3$ immediately below
the aforementioned occurrence of $E_2$ in order to prevent the edge from
exiting the box at point $3'$.  This procedure is only sustainable with
respect to Lemma \ref{secb3} if $n$
is odd, and in this case, the situation is as shown in Figure
\ref{fig2}, with the
extra condition that the vertical edge in the top left of the picture is
$L$-decorated.  There are two ways this picture can continue to the bottom
left consistently with Lemma \ref{secb3}:
 either the edge exits the box at point 
$1'$ without encountering further generators, or the edge encounters an 
occurrence of $E_1$.  The first situation cannot occur because it contradicts
Lemma \ref{secb4} and the hypothesis that $E \in \rdescent{\bu}$.  The second
situation cannot occur because it shows that $\bu$ is commutation equivalent
to a monomial of the form $\bv JIJ \bv'$, which contradicts the hypothesis
that $\bu$ be reduced.

We conclude that $\bu$ contains precisely
one occurrence of $E$.  By Lemma \ref{secb4} and the fact that $E \in
\ldescent{\bu} \cap \rdescent{\bu}$, this can
only happen if $\bu$ contains no occurrences of $E_1$.  

A similar argument shows that if $\bu'$ contains an occurrence of the 
generator $F$, then $\bu$ contains at most one occurrence of $F$, and it 
can only contain $F$ if it contains no occurrences of $E_{n-1}$.

It follows that at least one of the four situations must occur:
\begin{enumerate}
\item[(a)]{$\bu'$ contains $E$ and $\bu = E D_E$, 
where $D_E$ contains no occurrences of $E$ or $E_1$;}
\item[(b)]{$\bu'$ contains $F$ and $\bu = D_F F$, 
where $D_F$ contains no occurrences of 
$E_{n-1}$ or $F$;}
\item[(c)]{$\bu'$ contains neither $E$ nor $F$.}
\end{enumerate}

In cases (a) and (b), there is a corresponding factorization of $\bu'$, 
and the result claimed now follows from the faithfulness of the 
blob algebra as a diagram
calculus for the 
type-$B$ TL algebra \cite{blobcgm,marblobpres}.
For example, in case (a), we have
$\bu' = E D'_E$, and the fact that $E$ commutes with each generator in
each of $D_E$ and $D'_E$ implies that $D_E$ and $D'_E$ map to the same 
blob diagram, where the blob in this case is identified with $F$.

Suppose that we are in case (c), but that $\bu$ contains an occurrence
of $E$ or $F$.  Because the diagram $D$ corresponds to $\bu'$, it cannot 
have loops, so it must be the case that $\bu$ is either $L$-decorated at some
point, or $R$-decorated at some point.  This contradicts the hypotheses
on $\bu'$, using Lemma \ref{seca4}.  Since neither $\bu$ nor $\bu'$ contains
$E$ or $F$, the result follows by the faithfulness of the diagram calculus
for the Temperley--Lieb algebra 
\cite[\S6.4]{marbk}.
This completes the proof of (i).

We now prove (ii) in the case where $\bu' = IJ$; the case $\bu' = JI$
follows by a symmetrical argument.  
Thus, $\bu$ maps to the 
same diagram as $IJ$.  The fact that 
$\ldescent{\bu}$ is the set of generators in $I$ and $\rdescent{\bu}$ 
is the set of generators 
in $J$ means that $\bu$ cannot be left or right reducible.  
By Proposition \ref{secc3} (ii), this 
immediately means that $\bu = IJ$.
\end{proof}

\subsection{Proof of the theorem} 

\begin{lem}\label{secc4}
Let $\bu$ be a reduced monomial and let $D$ be the corresponding diagram.
Then $D$ avoids all the features on the left hand sides of Table
\ref{blobtab}.
Furthermore, $D$ contains at most one arc
with more than one decoration.
\end{lem}

\begin{proof}
The proof is by induction on the length of $\bu$.
If $\bu$ is a product of commuting generators, or $\bu = IJ$, or $\bu = JI$,
the assertions are easy to check, so we may assume that this is not the case.
(This covers the base case of the induction as a special case.)

By Proposition \ref{secc3}, $\bu$ must either be left reducible or right reducible.
We treat the case of left reducibility; the other follows by a symmetrical
argument.

By applying commutations to $\bu$ if necessary, we may now assume that 
$\bu = st \bv$, where $s$ and $t$ are noncommuting generators, and $t \not\in
\{E, F\}$.  By induction, we know that the reduced monomial $t \bv$ corresponds
to a diagram $D'$ with none of the forbidden features and at most one edge
with two decorations.

Suppose that $t = E_1$ and $s = E$.  By Lemma \ref{secb4}, points $1$ and $2$ of
$D'$ must be connected by an undecorated edge, and the effect of multiplying
by $E$ is simply to decorate this edge.  This does not introduce any forbidden
features, nor does it create an edge with two decorations, and this completes
the inductive step in this case.

The case where $t = E_{n-1}$ and $s = F$ is treated similarly to the above
case, so we may now assume that $s, t \not\in \{E, F\}$.  We must either have
$s = E_i$ and $t = E_{i+1}$, or vice versa.

Suppose that $s = E_i$ and $t = E_{i+1}$.  By
Lemma \ref{secb4}, this means that points $i+1$ and $i+2$ of $D'$ are connected 
by an undecorated edge.  The effect of multiplying by $s$ is then (a) to remove
this undecorated edge, then (b) to disconnect the edge emerging from point $i$ 
of $D'$ and reconnect it to point $i+2$, retaining its original decorated
status, then (c) to install an undecorated edge between points $i$ and
$i+1$.  This procedure does not create any forbidden features, nor does it
create a new edge with more than one decoration.

The case in which $s = E_{i+1}$ and $t = E_i$ is treated using a parallel
argument, and this completes the inductive step in all cases.
\end{proof}

\begin{lem}\label{secc5}
Let $\bu$ be a reduced monomial with corresponding diagram $D$.
\begin{enumerate}
\item[\rm (i)]{If points $1$ and $2$ (respectively, $1'$ and $2'$) 
are connected in $D$ by an edge decorated by $L$ but not $R$,
then $\bu$ is equal (as an algebra element) to a word of
the form $\bu' = E E_1 \bv$ (respectively, $\bu' = \bv E_1 E$).}
\item[\rm (ii)]{If points $n-1$ and $n$ (respectively, $(n-1)'$ and $n'$)
are connected in $D$ by an edge
decorated by $R$ but not $L$,
then $\bu$ is equal (as an algebra element) to a word of
the form $\bu' = F E_{n-1} \bv$ (respectively, $\bu' = \bv E_{n-1} F$).}
\end{enumerate}
\end{lem}

\begin{proof}
We first prove the part of (i) dealing with points $1$ and $2$.
By Lemma \ref{secb4}, we have $E \in \ldescent{\bu}$, so $\bu = E \bv'$.
Now $\bv'$ is also a reduced monomial, and by Lemma \ref{secc4}, $\bv'$
corresponds to a diagram $D'$ with no forbidden features.  Since multiplication
by $e$ does not change the underlying shape of a diagram (ignoring the
decorations), it must be the case that points $1$ and $2$ of $D'$ are 
connected by some kind of edge.  Since $D$ has no forbidden features and
the corresponding edge in $D$ has no $R$-decoration, the only way for this
to happen is if the edge connecting points $1$ and $2$ in $D'$ is undecorated.
By Lemma \ref{secb4}, this means that $\bv'$ is equal as an algebra element to
a monomial of the form $E_1 \bv$, and this completes the proof of (i) in
this case.

The other assertion of (i) and the assertions of (ii) follow by parallel 
arguments.
\end{proof}

\begin{proof}[Proof of Theorem \ref{thmmain}]
It is enough to prove the statement using the rescaling 
in which $\kl = \kr = 1$ (see section~\ref{sect:params} for details).

It is clear from the generators and relations that the reduced monomials are
a spanning set, and that the diagram algebra is a homomorphic image of the
abstractly defined algebra.  By Lemma \ref{secc4}, all reduced monomials map to
basis diagrams.  The only way the homomorphism could fail to be injective
is therefore for two reduced monomials $\bu$ and $\bu'$ to map to the same
diagram $D$, and yet to be distinct as algebra elements.  

It is therefore enough to prove that if $\bu$ and $\bu'$ are reduced monomials
mapping to the same diagram, then they are equal as algebra elements.
Without loss of generality, we assume 
that $\ell(\bu) \leq \ell(\bu')$ (where $\ell$ denotes length).

We proceed by induction on $\ell(\bu)$.  If $\ell(\bu) \leq 1$, or, more
generally, if $\bu$ is a product of commuting generators, then Lemma
\ref{secb5}
shows that $\bu = \bu'$.  Similarly, if $\bu = IJ$ or $\bu = JI$, then
$\bu = \bu'$, again by Lemma \ref{secb5}.  In particular, 
this deals with the base 
case of the induction.

By Proposition \ref{secc3}, we may now assume that $\bu$ is either 
left or right
reducible.  We treat the case of left reducibility, the other being similar.
By applying commutations if necessary, we may reduce to the case where
$\bu = st \bv$, $s$ and $t$ are noncommuting generators, and $t \not\in
\{E, F\}$.

Suppose that $s = E$, meaning that $t = E_1$.  In this case, points $1$ and
$2$ of $D$ are connected by an edge decorated by $L$ but not $R$.  By
Lemma \ref{secc5}~(i), this means that we have $\bu' = st \bv'$ as algebra
elements.
Since $\bu$ and $\bu'$ share a diagram, the (non-reduced) monomials
$t \bu$ and $t \bu'$ must also share a diagram.  Since $tst = \kl t = t$, 
the (reduced) monomials $t \bv$ and $t \bv'$ also map to the same diagram,
$D'$.  However, $t \bv$ is shorter than $\bu$, so by induction, 
$t \bv = t \bv'$, which in turn implies that $\bu = \bu'$.

Suppose that $s = F$, meaning that $t = E_{n-1}$.  An argument similar to
the above, using Lemma \ref{secc5}~(ii), establishes that $\bu = \bu'$ in this
case too.

We are left with the case where $s = E_i$ and either $t = E_{i+1}$ or
$t = E_{i-1}$ (where $t \not\in \{E, F\}$).  We will treat the case where
$t = E_{i+1}$; the other case follows similarly.  In this case, we have
$tst = t$, and so $t \bu = tst \bv = t \bv$.  It is not necessarily true that
$t \bu'$ is a reduced monomial, but it maps to the same diagram as $t \bv$,
which is reduced.  After applying algebra relations to 
$t \bu'$, we may transform it into a scalar multiple of a reduced monomial,
$\br$.  Since reduced monomials map to basis diagrams (Lemma \ref{secc4}),
the scalar involved must be $1$.  Now the reduced monomials $t \bv$ and
$\br$ map to the same basis diagram, and $t \bv$ is shorter than $\bu$,
so by induction, we have $t \bv = \br$ as algebra elements.

Since $s \in \ldescent{\bu}$, we have $s \in \ldescent{\bu'}$ by Lemma
\ref{secb4}, so that $\bu' = s \bv''$ for some reduced monomial $\bv''$.
Since $sts = s$, we have $s (t \bu') = \bu'$.  We have shown that
$t \bu' = \br = t \bv$, so we have $$
\bu' = s t \bu' = st \bv = \bu
,$$ which completes the proof.
\end{proof}


\section{A cellular basis}
\newcommand{\ASTLA}{affine symmetric Temperley--Lieb algebra}

In what follows we shall make repeated use of the cellular basis for
the symplectic blob algebra (or equivalently the \ASTLA). 
In this section we review the construction of this cellular basis. 


Let $d \in B_{2n}^\phi$ and define $\#(d)$ to be the number of
propagating lines in $d$. We extend this notation to scalar multiples
of $d$.

Suppose $\#(d) \ge 2$.  Then there is a unique pair of propagating
lines that can be simultaneously deformed to touch the 0-reflection
line:
$$ \epsfbox{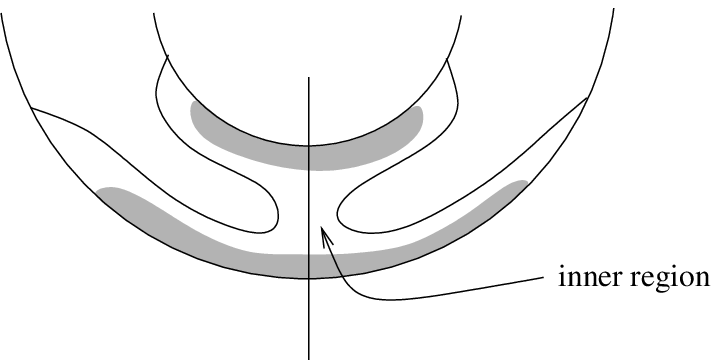}$$
the ``closest'' propagating lines to the 0-reflection line. This
defines a unique \emph{inner region} that will be black or white when
the diagram is coloured.

We define:
$$c(d) = \left\{ \begin{array}{ll} b &\mbox{if the inner region is
      black}\\
    w &\mbox{if the inner region is white}\end{array}\right.$$

\begin{lem}[{\cite[lemma 8.2.1]{gensymp}}]
For all $d, d' \in B_{2n}^{\phi}$ we have:
\begin{enumerate}
\item[(i)] $\#(dd') \le \#(d)$
\item[(ii)] If $\#(dd') = \#(d)$ and $\#(d) \ne 0$ then $c(dd') = c(d)$.
\end{enumerate}
\end{lem}

Inspired by this we define for $0 
< l \le n$:
$$B_{2n}^{\phi}[l] = \{ d \in B_{2n}^{\phi}\mid \#(d)=2l \mbox{ and }
c(d) =b \}$$
$$B_{2n}^{\phi}[-l] = \{ d \in B_{2n}^{\phi}\mid \#(d)=2l \mbox{ and }
c(d) =w \}$$
and
$$B_{2n}^{\phi}[0] = \{ d \in B_{2n}^{\phi}\mid \#(d)=0\}.$$

%
%
%
%
%

We define for $-n \le l \le n-1$
$$B_{2n}^{\phi}(l) = B_{2n}^{\phi}[l] \cup \bigcup_{-|l| < a < |l|}
B_{2n}^{\phi}[a]$$
and $I_{2n}^{\phi}(l)$ to be the ideal of $b_{2n}^\phi$ generated by
$B_{2n}^{\phi}[l]$

\begin{prop}[{\cite[proposition 8.2.2]{gensymp}}]
The ideal $I_{2n}^{\phi}(l)$  has basis
$B_{2n}^{\phi}(l)$.
\end{prop}

We have set inclusions:
$$
\xymatrix{
          &  B^{\phi}_{2n}(n-1) \ar @{_(->}[dl] <1ex>  
                 &\cdots  \ar@{_{(}->}[l] \ar@{_(->}[ddl] <1ex> 
&  B^{\phi}_{2n}(2) \ar@{_{(}->}[l] \ar@{_(->}[ddl] <1ex> 
 &  B^{\phi}_{2n}(1) \ar@{_{(}->}[l] \ar@{_(->}[ddl] <1ex>  
\\
B^{\phi}_{2n}(-n) & & & 
          & &  B^{\phi}_{2n}(0) \ar@{_{(}->}[ul]
                                       \ar@{_{(}->}[dl] <1ex>  
\\
          &  B^{\phi}_{2n}(-(n-1)) \ar@{_{(}->}[ul] 
               & \cdots  \ar@{_{(}->}[l] \ar@{_(->}[uul]   
 &  B^{\phi}_{2n}(-2) \ar@{_{(}->}[l] \ar@{_(->}[uul] <1ex> 
  &  B^{\phi}_{2n}(-1)) \ar@{_{(}->}[l] \ar@{_{(}->}[uul] }
$$
and this passes to a subideal structure (\cite[proposition 8.2.2]{gensymp}).
We can use the ideals $I_{2n}^\phi(l)$ to define a cellular structure on
$b_{2n}^{\phi}$.
We define $S_l=I_{2n}^{\phi}(-l)$ and $T_l= S_l + I_{2n}^\phi(l)$.
We have a chain of ideals:
$$b_{2n}^{\phi} = S_n \supseteq T_{n-1} \supseteq
S_{n-1} \supseteq T_{n-2} \supseteq \cdots \supseteq T_1 \supseteq
S_1 \supseteq S_0.
$$

\begin{thm}[{\cite[theorem 8.2.8]{gensymp}}]
The above chain of ideals is a cellular chain for $b_{2n}^{\phi}$ and
the set $B^\phi_{2n}$ is a cellular basis for $b_{2n}^{\phi}$.
\end{thm}

We define for $l \in \{-n, -n+1, \ldots, n-1 \}$ and $d \in
B_{2n}^{\phi}[l]$:
$$\sS_{2n}^d(l) := \frac{b_{2n}^{\phi} d + T_{|l| -1}}{T_{|l| -1}}.$$
These are the cell modules for $b_{2n}^\phi$ 
(which have the flavour of Specht modules).

\begin{thm}[{\cite[theorem 8.2.9]{gensymp}}]
If $\delta$, $\delta_L$, $\delta_R$, $\kappa_L$, $\kappa_R$, and 
$\kappa_{LR}$ are
all units then $b_{2n}^\phi$ is quasi-hereditary with heredity chain
as above.
\end{thm}

Thus in the case where all the parameters are units then the cell
modules defined above are standard modules which depend only on $l$
(and not on $d$) and there are exactly $2n$
simple modules for $b_{2n}^{\phi}$.

We may label the simple modules with the following poset:
$$
\begin{array}{c}
\xymatrix@R=8pt@C=4pt{
&0 \ar@{-}[dl] \ar@{-}[dr]& \\
1 \ar@{-}[d]\ar@{-}[drr] & &-1 \ar@{-}[d]\ar@{-}[dll]\\
2 \ar@{-}[d]\ar@{-}[drr] & &-2 \ar@{-}[d]\ar@{-}[dll]\\
{\genfrac{}{}{0pt}{}{\vdots}{\vdots}} \ar@{-}[d]\ar@{-}[drr] & 
&{\genfrac{}{}{0pt}{}{\vdots}{\vdots}} \ar@{-}[d]\ar@{-}[dll]\\
n-2 \ar@{-}[d]\ar@{-}[drr] & &-n+2 \ar@{-}[d]\ar@{-}[dll]\\
n-1 \ar@{-}[dr] & &-n+1 \ar@{-}[dl]\\
&-n& 
}
\end{array}
$$
where $0$ is the maximal element and $-n$ is the minimal element.
When the algebra is quasi-hereditary with the above poset, this has
important consquences for the representation theory of the algebra.

%
%
%

\begin{rem}
The parameter $\kk$ can only appear
in
expressions where there are diagrams with no propagating lines.
Thus, it has no effect on the action of the algebra on standard
modules whose basis has at least one propagating line.
Thus, to study standard modules with non-zero label we could safely ignore the
value of $\kk$ and not affect the representation theory.
We have several successive quotients:
$$ H(\tilde{C}) \to \twoBTL \to b_n^x \tilde{\to}  
b_{2n}^\phi \to \frac{b_{2n}^\phi}{I^\phi_{2n}(0)}$$
where $H(\tilde{C})$ is the Hecke algebra of type $\tilde{C_n}$,
$I^\phi_{2n}(0)$ is the ideal in $b_{2n}^\phi$ generated by elements
with no propagating lines, and $\twoBTL$ is defined in section 
\ref{sect:pres}.
The algebra $\frac{b_{2n}^\phi}{I^\phi_{2n}(0)}$ is independent of $\kk$.
\end{rem}

The $b_{2n}^\phi$ version of 
proposition \ref{bxeiso} is the following:
\begin{prop}[{\cite[proposition 8.1.1]{gensymp}}]\label{eiso}
Suppose $\delta_L \ne 0$, and set $e' = \frac{e}{\delta_L}$.  Then we have
$$e' b_{2n}^{\phi} (\delta, \delta_L, \delta_R, \kappa_L, \kappa_R,
\kappa_{LR}) e' 
\cong b_{2(n-1)}^{\phi}(\delta, \kappa_L, \delta_R, \delta_L,
\kappa_R, \kappa_{LR}).
$$
\end{prop}

We also have a ``dual'' version of proposition \ref{eiso}
(the $b_{2n}^\phi$ version of 
proposition \ref{bxfiso}).
\begin{prop}\label{fiso}
Suppose $\delta_R \ne 0$, and set $f' = \frac{f}{\delta_R}$. Then we have
$$f' b_{2n}^{\phi} (\delta, \delta_L, \delta_R, \kappa_L, \kappa_R,
\kappa_{LR}) f' 
\cong b_{2(n-1)}^{\phi}(\delta, \delta_L, \kappa_R, \kappa_L,
\delta_R, \kappa_{LR}).
$$
\end{prop}

We can thus use localisation and globalisation functors to inductively
prove many properties for $b_{2n}^\phi$. (Our choice of relations was
partly designed to make this proposition work, but still obtain an
accessible quotient of the Hecke algebra of type $\tilde{C}$.)
These functors will be introduced later in section \ref{except}.




\section{Good parametrisations}\label{sect:params}

The symplectic blob algebra is a quotient of the Hecke algebra of type
$\tilde{C}$, which has three parameters
(itself a quotient of the group algebra of a Coxeter--Artin system). 
Yet our definition has six.
One of the parameters, $\kk$, is used to make the quotient `smaller' 
than the Hecke algebra (indeed finite dimensional),
just as the three Hecke parameters make the Hecke algebra smaller than
the braid group. 
Provided  the appropriate parameters are units, two of
the other parameters can  be scaled away, 
as we shall see, leaving  a parameter set 
corresponding to that of the Hecke algebra. 

However, our aim here is to determine the exceptional representation
theory of the symplectic blob (and hence part of the Hecke/braid 
representation theory). Therefore we are interested in working over
rings from which we can base change to the exceptional cases. 
Clearly this requires some up-front knowledge of what the exceptional
cases are. This bootstrap problem can be solved in significant part
by looking at Gram matrices for standard modules, exactly
as in \cite{martsaleur,marwood}.
We give details in Section~\ref{ss:gram}, 
but the good blob algebra parametrisation 
\cite{marryom}
leads us
to consider four parameters $q,l,r$ and $\kk$, determining
$\dd = [2]$ 
$\dl = [l]$;
$\dr= [r]$;
$\kl=[l-1]$;
$\kr=[r-1]$.
Here the notation $[n]$ is used to denote
the $q$-number $\frac{q^n -q^{-n}}{q-q^{-1}}$, with $q \in k$.
(If $k = \C$ then it is usual to take $q = \exp(i\pi /m)$, thus defining
an equivalent parameter $m$.)

The new parameters $q$, $l$, $r$ may not even be real, if one wishes to work
in the complex setting. But if $l$ and $r$ are integral then we can work in
the ring of Laurent polynomials (and the integral cases are the most
singular, as flagged by the Gram determinants).



\newcommand{\eql}[1]{\begin{equation}\label{#1}}
\newcommand{\eq}{\end{equation}}

\subsection{Gram matrices} \label{ss:gram}

\begin{figure}
\caption{\label{ppm01}}
\[
\epsfysize=60mm \epsfbox{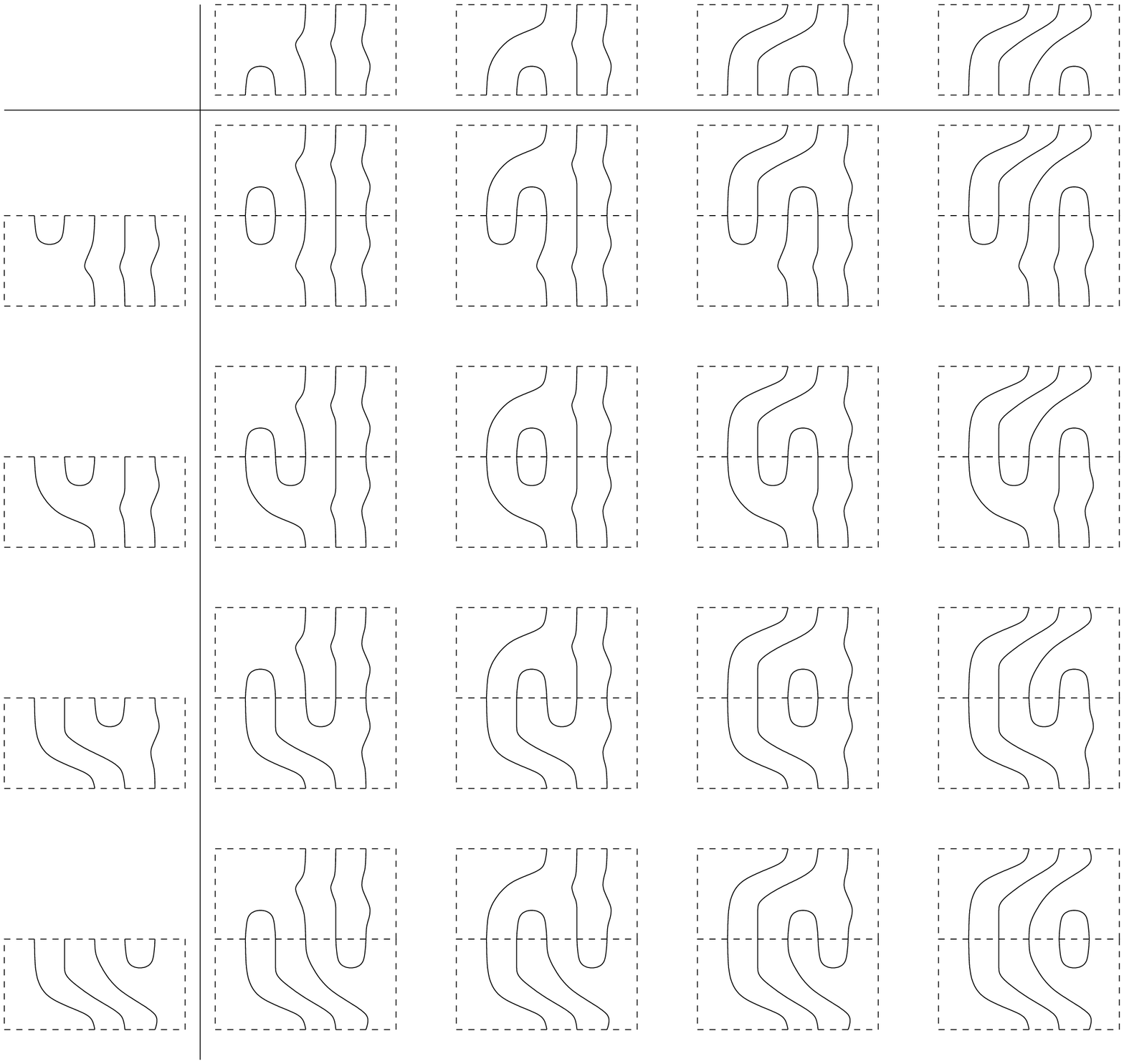}
\]
\end{figure}

Given an algebra $A$ with an involutive antiautomorphism, and an
$A$-module with a simple head and a natural inner product, then a
condition for semisimplicity of a specialisation 
is that this inner product be non-degenerate. 
Conversely, for non-semisimple cases, we look for
conditions under which the Gram determinant vanishes. 
Depending on the number and incarnation of parameters, this vanishing
may appear to describe a complicated variety. 
However experience shows
that for a {\em good} choice of incarnation the non-semisimplicity condition
can often be stated simply. 
For an initial illustrative example, 
consider  Temperley--Lieb cell modules in the `upper half-diagram' bases, 
such as the following $n=5, l=3$ case. 
The basis (acted on by diagrams from above, ignoring irrelevant arcs
below) may be written
\[
H_T^3(5,3): \qquad
\epsfysize=8mm \epsfbox{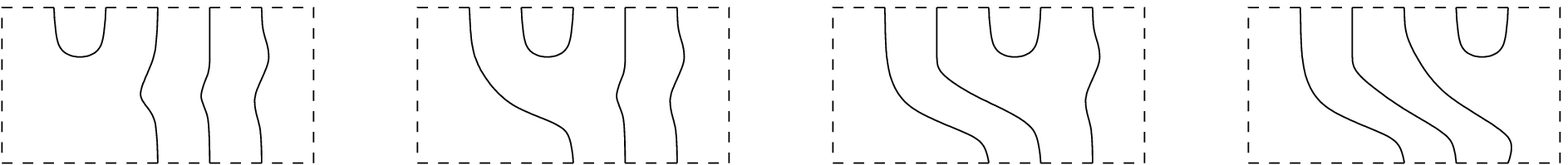}
\]
whereupon the inner product is computed using the array in 
Figure~\ref{ppm01} (via the usual Temperley--Lieb diagram inversion antiautomorphism).
This gives immediately the Gram matrix
\newcommand{\gM}{M^{TL}} 
\newcommand{\mat}[1]{\left( \begin{matrix}}
\newcommand{\tam}{\end{matrix} \right)}
\[
\gM_5(3) = \mat{ccccc}
\delta & 1 & 0 & 0 \\
1 & \delta & 1 & 0 \\
0 & 1 & \delta & 1 \\
0 & 0 & 1 & \delta
\tam
\]
(note the slight difference in the way diagrams act on the basis, 
cf. usual diagram multiplication). 
The generalisation to $\gM_n(n-2)$ will be obvious. 
Evaluation of the determinant is also straightforward. 
\newcommand{\oplo}{\oplus_1^1}
Indeed, more generally still, 
writing $M' \oplo M''$ for the almost block diagonal matrix 
\[ 
M' \oplo M''  \; := \;
\left(
\begin{array}{ccc|ccc}
\cline{1-3}
 & & \\
 &M'& & &0 \\
 & & & 1 \\ \cline{1-6}
&&1 &&&\\
&0 && & M'' & \\
&&&&& \\
\cline{4-6}
\end{array}
\right) \]
and $\mu_n(M) = M \oplo (\delta) \oplo (\delta) ...\oplo (\delta)$ 
($n+1$ terms) for any initial matrix $M$ so, for example, that 
$\gM_n(n-2)  = \mu_{n-2}((\delta))$, 
we have
\begin{equation} \label{eq:p2}
\det(\mu_n(M))= \delta \; \det(\mu_{n-1}(M)) - \det(\mu_{n-2}(M))
\qquad (n>0)
\end{equation}
where $\mu_{-1}(M) = M^{dd}$ (the matrix $M$ with the last row and
column removed). 


As is well known the recurrence 
\begin{equation} \label{eq:p3}
M(n)=[2]M(n-1) -M(n-2)
\end{equation}
is solved by $M(n) = \alpha [s+n]$ for any constants $s, \alpha$. 
(Two pieces of initial data, such as $M(0), M(1)$, fix them via 
$M(0)=\alpha [s]$, $M(1)= \alpha [s+1]$.)

Comparing (\ref{eq:p3}) with  (\ref{eq:p2})
then leads us to the parametrisation $\delta=[2]$, which 
 makes $ \det(\gM_n(n-2)) = [n]$. The vanishing of this form is
well understood, requiring $q$ to be a root of 1. 
(Although this only gives one Gram determinant per algebra, abstract
representation theory tells us that it gives a complete picture in the
Temperley--Lieb case.) 

\begin{figure}
\caption{\label{ppm1}}
\[
\epsfysize=80mm \epsfbox{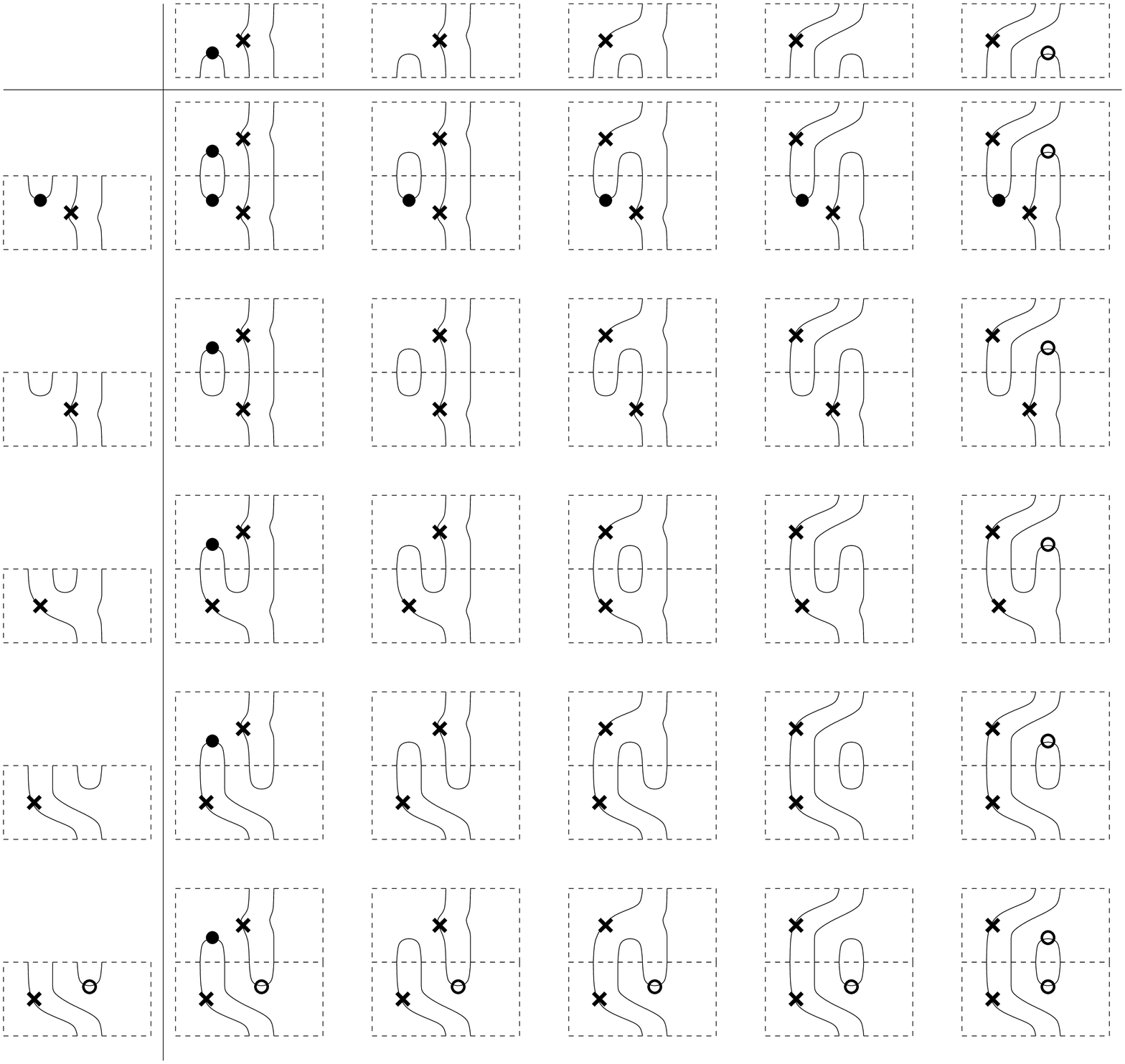}
\]
\end{figure}

A similar analysis proceeds for the ordinary and symplectic blob
algebras. For example a basis for one of our symplectic cell modules
is the   
left-most (labelling) column of the array in Figure~\ref{ppm1}. 
In fact this picture simultaneously encodes three of our $n=4$ cell modules
$\sS_{2n}^{}(l)$ (those with label $l=-2,1,-1$),
depending on how blobs are understood to act on the two propagating
lines (the leftmost of which is marked with a $\times$ to flag this choice). 
Indeed, by omitting the last row and column we get the corresponding
array for (two cell modules of) the ordinary blob algebra. 

In the blob case, 
choosing a parametrisation in which $e$ is idempotent, 
for arithmetic simplicity, 
we thus have  Gram matrices $M^b_n(n-2) = \mu_{n}(B_+)$, 
and  $M^b_n(-(n-2)) = \mu_{n-1}(B_-)$, 
where
\[
B_+ = \mat{cc} \kl & \kl \\ \kl & \delta \tam,
\qquad
B_- = \mat{ccc} \kl & \kl & 1 \\  \kl & \delta &1 \\ 1&1&\delta \tam
\]
In other words we have the same bulk recurrence as for Temperley--Lieb, 
but more interesting initial conditions. 
The idea is to try to parametrise $\kl$ so that the initial conditions
conform to the natural form $M(n)=\alpha [s+n]$, 
for some choice of $\alpha,s$, in each case. 
We have 
\[
(B_+): \;\;
\kl = \alpha [s+1], \qquad \kl([2]-\kl)  = \alpha [s+2],
\]
Eliminating $\alpha$, one sees
that the parametrisation $\kl=\frac{[l]}{[l+1]}$ 
is indicated 
(or equivalently $\kl=\frac{[l]}{[l-1]}$
by exchanging $l$ and $-l$).
We have chosen the symbol $l$ for the $s$-parameter in $\kl$ for obvious reasons. 
The same parametrisation works for $B_-$. 
This gives
$\det(M^b_n(n-2))= \frac{[l]}{[l+1]^2} [n+l]$ and 
$\det(M^b_n(-(n-2)))= \frac{[l+2]}{[l+1]^2} [2-n+l]$. 
Once again this is in a convenient form to simply characterise the
singular cases for all $n$.
The result (as is well known \cite{marryom}) 
is a generalised form of the kind of
alcove geometry that occurs in Lie theory
(or more precisely in quantum group representation theory when $q$ is
an $m$-th root of 1). 
The blocks of the algebra
are described by orbits of an affine reflection group, with the
separation of affine walls determined by $m$, and the `$\rho$-shift' of the
origin determined by $l$. 

Returning finally to the symplectic case, 
choosing a normalisation in which $e$ and $f$ are idempotent, 
the most complicated of the three cell choices gives 
(in the obvious generalisation to any $n$, and $n-2$ propagating
lines) the $(n+1) \times (n+1)$ matrix:
\eql{sblob gram1}
M'_{}(n,\kl,\kr ) = \mat{cccccccccccc}
\kl & \kl & 1&0 &0&0& \cdots &0 \\
\kl & [2] & 1 & 0  &0&0& \cdots &0\\
1 & 1 & [2] & 1 & 0 &0& \cdots &0 \\
0 & 0 & 1 & [2] & 1 & 0 & \cdots &0 \\
\vdots&&&& \ddots & & &\vdots \\
0 & \cdots & 0 & 0 & 1 & [2] &1& 0 \\
0 & \cdots & 0 & 0 & 0 & 1 & [2]& \kr \\
0 & \cdots & 0 & 0 & 0 & 0 & \kr  & \kr 
\tam
\eq
This is clearly just another variant of the Temperley--Lieb Gram matrix, with yet more
interesting `boundary conditions'. 
Laplace expanding we have:  
\[
\det(M'_{}(n,\kl,\kr )) = \kr \frac{[l+2]}{[l+1]^2} ([2-n+l] - \kr [3-n+l]).
\]
Again one seeks to parametrise $\kr$ by  a parameter $r$ 
such that, for every $n$, this determinant becomes a simple product of
quantum numbers, quantum-integral when $r$ is integral. 
Clearly $\kr = \frac{[r]}{[r+1]}$ does this 
(the precise form is chosen for symmetry with $\kl$). One obtains:
\[
\det(M'_{}(n,\kl,\kr )) 
= \frac{[r][l+2]}{[r+1][l+1]} 
      \frac{[l-(r+n-2)]}{[r+1][l+1]}
\]
(and similarly for the others
--- we give complete results for these cases and more,
and demonstrate that our normalisation choice was without loss of generality,
in Section~\ref{sect:gramdets}). 
Once again then, the most singular cases are when $q$ is a root of 1
and $l$ and $r$ are integral.
It is intriguing to speculate on the corresponding alcove geometry
(cf. the ordinary blob case above) --- see later. 
However there is a significant difference in this symplectic case, in that
here abstract representation theory does not allow us simply to reconstruct all
other Gram determinants from this subset. 
It is 
for this reason that we return to report some less straighforward Gram
determinant results in Section~\ref{sect:gramdets}.


\subsection{Four rescalings}
We tabulate the parameter shifts associated to four choices of generator
rescalings, such that only four free parameters remain in each case: 
\newcommand{\mo}{\mapsto \;\;\;}
\newcommand{\ignore}[1]{}
\[
\begin{array}{l|rr|rrrrrrl}
    &\multicolumn{2}{c}{\mbox{gen. shifts} }
  &\multicolumn{6}{c}{\mbox{parameter shifts} }
\\ \hline 
\mbox{way} & e\mo   &   f\mo   &  \dd\mo &  \dl\mo &  \dr\mo & \kl\mo&\kr\mo&\kk\mo 
\\ \hline 
1 &\frac{e}{\dl}&\frac{f}{\dr}& \dd
  &1&1&\frac{\kl}{\dl}&\frac{\kr}{\dr}&\kk/(\dl\dr)
\\
2 &\frac{e}{\kl}&\frac{f}{\kr}&\dd&\frac{\dl}{\kl}&\frac{\dr}{\kr}&1&1&\kk/(\kl\kr)
\\
3 &\frac{e}{\kl}&\frac{f}{\dr}& \dd &\frac{\dl}{\kl}&1&1&\frac{\kr}{\dr}&
 \kk/(\dr\kl)
\\
4 &\frac{e}{\dl}&\frac{f}{\kr}&\dd& 1 &
\frac{\dr}{\kr}&\frac{\kl}{\dl}& 1 &
 \kk/(\dl\kr)
\\ \hline 
\mbox{blob}&&&[2]&[l]&[r]&[l-1]&[r-1]&\kk
\\
\mbox{DN}& & & [2] &\frac{[\omega_1]}{[\omega_1 +1]}
  &\frac{[\omega_2]}{[\omega_2 +1]} & 1&1&b
\\
\mbox{GMP}&&&[2]&[w_1]&[w_2]&[w_1+1]&[w_2+1]&\kk
\end{array}
\]

This is easy to check. 
For instance, we know that the symplectic blob has generators that satisfy
the relations given in the presentation in section \ref{sect:pres},
and since the parameter
change is given by rescaling the generators, it is easy to check what
the new values of the parameters should be.

Note, however that only using four parameters obscures the swapping of 
the parameters induced by the globalisation and localisation 
functors $G$, $G'$, $F$ and $F'$ which will be introduced later.

\subsection{Two more ways to reparametrise}\label{ourpara}
The `good' parametrisation is recalled in the row labelled {\em blob} in our
table. 
Scaled 
 in ``Way 2'', this is similar to the parametrisation used by De Gier and
Nichols in \cite{degiernichols} (row DN). 
\ignore{{
They have:
\begin{align*}
\dd&=[2] & \dl &=\frac{[\omega_1]}{[\omega_1 +1]}\\
\kl &=1 & \dr &=\frac{[\omega_2]}{[\omega_2 +1]}\\
\kr &=1 & \kk &=b
\end{align*}
}}
They in effect have four parameters, $q$, $\omega_1$, $\omega_2$
and $b$,
but significantly they reparametrise $b$ in terms of box numbers
and $\theta$.
We will use the parametrisation in row GMP
\ignore{{
\begin{align*}
\dd&=[2] & \dl &=[w_1]\\
\kl &=[w_1+1] & \dr &=[w_2]\\
\kr &=[w_2+1] &  &
\end{align*}
and
}}
with
$$
\kk = \left\{ \begin{array}{ll}
      {\displaystyle{\left[{\frac{w_1+w_2+\theta +1}{2}}\right]
      \left[{\frac{w_1+w_2-\theta +1}{2}}\right]}} &\mbox{if $n$
        even}\\
         & \\
      {\displaystyle{-\left[{\frac{w_1-w_2+\theta }{2}}\right]
      \left[{\frac{w_1-w_2-\theta }{2}}\right]}} &\mbox{if $n$
        odd.}
\end{array}\right.
$$
which is connected to the parametrisation of De Gier and
Nichols via ``Way 2''. 
We use GMP 
to avoid too
many terms in the denominator for the Gram determinants calculated
later.

%




\section{On exceptional cases}\label{except}
Throughout this section we will assume that all the parameters are
units, so that $b_{2n}^{\phi}$ is quasihereditary. 

\subsection{Globalisation and localisation functors}
We now define various globalisation and localisation functors. Their
construction is very general and further information about their
properties may be found in \cite[section 2]{gensymp}, \cite{green} or
\cite{marryom}.

Recall that $e' = \frac{e}{\dl}$ and $f' = \frac{f}{\dr}$.
Let
\begin{align*}
G:e' b_{2n}^{\phi} e'\mbox{-}\Mod &\to b_{2n}^{\phi}\mbox{-}\Mod\\
M &\mapsto b_{2n}^{\phi} e' \otimes_{e'b_{2n}^{\phi} e'}M  
\end{align*}
be the globalisation functor with respect to $e$ and
\begin{align*}
G':f' b_{2n}^{\phi} f'\mbox{-}\Mod &\to b_{2n}^{\phi} \mbox{-}\Mod\\
M &\mapsto b_{2n}^{\phi} f' \otimes_{f'b_{2n}^{\phi} f'}M  
\end{align*}
be the globalisation functor with respect to $f$.
The functors $G$ and $G'$ are both right exact and $G\circ G' = G'
\circ G$.

We also have localisation functors:
\begin{align*}
F:b_{2n}^{\phi}\mbox{-}\Mod &\to e' b_{2n}^{\phi} e'\mbox{-}\Mod \\
M &\mapsto  e'M  
\end{align*}
\begin{align*}
F':b_{2n}^{\phi}\mbox{-}\Mod &\to f' b_{2n}^{\phi} f'\mbox{-}\Mod \\
M &\mapsto  f'M.
\end{align*}
The functors $F$ and $F'$ are both exact. They also map simple modules
to simple modules (or zero) \cite{green}.
We have $F\circ G = \id$ and $F'\circ G' =\id$.

We will abuse notation slightly and identify a module for the algebra 
$e' b_{2n}^{\phi} e'$ with its image via the isomorphism of
Proposition~\ref{eiso}
as a module for $b_{2n-2}$ and similarly for $f' b_{2n}^{\phi} f'$. The
functors applied will make it clear which parameter choice we are
making.

Standards behave well with respect to globalisation \cite[proposition
8.2.10]{gensymp}:
$$G\sS_{2n-2}(l) = \sS_{2n}(-l)$$ 
$$G'\sS_{2n-2}(l) = \sS_{2n}(l).$$ 
If we have a standard module $\sS_{2n}(l)$ that is not simple then
there is a simple module $L_{2n}(m)$ in the socle of $\sS_{2n}(l)$
with $m < l$, because $\sS_{2n}(l)$ is a standard module (and using
the definition of standard modules for quasi-hereditary algebras).
Since  $L_{2n}(m)$ is the head of $\sS_{2n}(m)$ 
there thus exists an $m < l$ such that there is a non-zero map
$$\sS_{2n-2}(m) \stackrel{\psi}\to \sS_{2n-2}(l).$$
Globalising then gives us:
$$\sS_{2n}(-m) \stackrel{G \psi}\to \sS_{2n}(-l)$$
$$\sS_{2n}(m) \stackrel{G' \psi}\to \sS_{2n}(l)$$
with $G\psi$ and $G'\psi$ both non-zero.
Note that the non-zero
map is non-zero on the simple head of the standard module, and hence the
head, $L_{2n}(m)$, must be a composition factor of the image in
$S_{2n}(\pm l)$. As $m <l$,
this factor is not equal to $L_{2n}(l)$ or $L_{2n}(-l)$, which implies
that
$\sS_{2n}(l)$ and $\sS_{2n}(-l)$ are also not simple. 
Thus parameter choices that give non-simple standards propagate.
(We must take care with the parameter swapping effect of $G$ and $G'$,
however.)

We set $L_{2n}(l)$ to be the irreducible head of the standard module
$\sS_{2n}(l)$.
We have the following proposition.
\begin{prop}\label{fonlands}
We have $$
FL_{2n}(l) = \begin{cases}
              0  &\mbox{if $l=-n$ or $l=-n+1$,}\\
              L_{2n-2}(-l) &\mbox{otherwise;}
             \end{cases}$$ 
$$
F'L_{2n}(l) = \begin{cases}
              0  &\mbox{if $l=-n$ or $l=n-1$,}\\
              L_{2n-2}(l) &\mbox{otherwise;}
             \end{cases}$$ 
$$F\sS_{2n}(l) = \sS_{2n-2}(-l) \quad\mbox{for $l \ne -n$ and $l \ne -n+1$};$$ 
$$F'\sS_{2n}(l) = \sS_{2n-2}(l)\quad\mbox{for $l \ne -n$ and $l \ne n-1$}.$$ 
\end{prop}
\begin{proof}
As $e$ or $f$ may be taken as part of a heredity chain
we obtain 
using \cite[proposition 3]{marryom} or \cite[appendix A1]{donkbk} the 
result for standard modules.
We may use 
\cite[proposition 2.0.1]{gensymp} and the above result for the
globalisation functor to determine which simple modules $F$ or $F'$
maps to zero. 
\end{proof}


\subsection{Some simple standard modules}\label{subsect:onedim}
By counting the number of diagrams in $B_{2n}[l]$ for $l=-n$, $-n+1\ne0$, 
$n-1\ne0$ (and $n-2\ne 0$) and 
considering the action of $b_{2n}^\phi$ it is clear that
\begin{lem}\label{standardsimple}
 Suppose that $n\ge2$, then 
$\sS_{2n}(-n)$, $\sS_{2n}(-n+1)$ and $\sS_{2n}(n-1)$ are all one-dimensional
and hence are irreducible. 
If further $n\ge 3$ then $\sS_{2n}(n-2)$ is one-dimensional
and hence irreducible. 
\end{lem}
When $b_{2n}^\phi$ is quasi-hereditary
and if $n \ge0 $ then
$\sS_{2n}(-n)=L_{2n}(-n)$ is the trivial module, i.e. the one-dimensional
module where all the generators of the algebra act as zero.
If $n \ge2 $ then
the module $\sS_{2n}(-n+1)=L_{2n}(-n+1)$ is the one-dimensional
module where $f$  acts as multiplication by
$\delta_R$ and
all the other generators of the algebra act as zero.
Similarly,
if $n \ge2 $ then
the module $\sS_{2n}(n-1)=L_{2n}(n-1)$ is the one-dimensional
module where $e$  acts as multiplication by
$\delta_L$ and
all the other generators of the algebra act as zero.

We may similarly note that 
if $n \ge3 $ then
the module $\sS_{2n}(n-2)=L_{2n}(n-2)$ is the one-dimensional
module where $e$  acts as multiplication by
$\delta_L$, $f$  acts as multiplication by
$\delta_R$ and
all the other generators of the algebra act as zero.

We use $\Ext^i(-,-)$ to denote the right derived functors of
$\Hom(-,-)$ which may be defined in the usual way in 
$\Mod b_{2n}^\phi$ as $b_{2n}^\phi$ is quasi-hereditary (and so there are
enough projectives and injectives).

We also use $[M:L]$ to denote the multiplicity of $L$, a simple
module as a composition factor of a (finite dimensional) module $M$.

We also note that the algebra $b_{2n}^\phi$ has a simple-preserving duality ---
namely the one induced by the algebra anti-automorphism that turns diagrams
upside down. 

We refer the reader to \cite[appendix A]{donkbk} or \cite{dlabring1}
for the definition and general properties of quasi-hereditary
algebras.

\begin{lem}\label{lem:noext}
If $b<a$ and  $[\sS_{2n}(a):L_{2n}(b)] =0$, then 
$$\Ext^1(L_{2n}(a), L_{2n}(b)) = 
\Ext^1(L_{2n}(b), L_{2n}(a)) =0.$$
\end{lem}
\begin{proof}
Assume for a contradiction that there is a non-split extension of
$L_{2n}(b)$ by $L_{2n}(a)$  for $a>b$ and  $[\sS_{2n}(a):L_{2n}(b)] =0$.
Now note that $\Ext^1(L_{2n}(a), L_{2n}(b))\cong \Ext^1(L_{2n}(b),
L_{2n}(a))$ as $b_{2n}^\phi$ has a simple-preserving duality. 
Thus we take $E$ to be the non-split extension defined by the
short exact sequence
$$ 0 \to L_{2n}(b) \to E \to L_{2n}(a) \to 0.$$
Now as $\sS_{2n}(a)$ is a standard module, it is the largest
$b_{2n}^\phi$ module with simple head $L_{2n}(a)$ and all other
composition factors having labels less than $a$.
Thus there must be a surjection $\sS_{2n}(a) \to E$.
This implies that $[\sS_{2n}(a), L_{2n}(b)] \ne 0$, the desired
contradiction.
\end{proof}

Thus, for $n\ge3$, by quasi-heredity, and the previous lemma
 there are no non-split 
extensions between the modules
$L_{2n}(-n)$,
$L_{2n}(n-1)$,
$L_{2n}(-n+1)$ and
$L_{2n}(n-2)$.

As  $b_{2n}^\phi$ has a simple-preserving duality, 
when we are coarsening the quasi-hereditary order we
need only consider the composition factors of the standard modules. 
In other words, if two adjacent labels $a > b$ in our original poset satisfy
$[\sS_{2n}(a), L_{2n}(b)] =0$, then the relation $a >b$ can be removed
from the poset and these labels do not
need to be comparable in a poset giving a quasi-hereditary order.
This procedure is discussed in greater detail after Lemma 1.1.1 in
\cite{erdpar}.
%

\subsection{Some composition multiplicities and coarsening the labelling poset}
\label{fourthree}

\begin{prop}\label{notcompfactor} 
Suppose $l \ge2 $ and $n \ge2$.  Then we have
$$[\sS_{2n}(l-1):L_{2n}(\pm l)] = 0$$
and 
$$[\sS_{2n}(-l+1):L_{2n}(\pm l)] = 0.$$
Suppose further that $l \ge3 $ and $n \ge 3$.  Then we have
$$[\sS_{2n}(l-2):L_{2n}(-l)] = 0,$$
and
$$[\sS_{2n}(-l+2):L_{2n}(l)] = 0.$$  (We interpret $[M : L_{2n}(n)]$ to
mean zero.)
\end{prop}
\begin{proof}
%
%
%

We first prove that $[\sS_{2n}(\pm(l-1)):L_{2n}(-l)] = 0$ by induction
on $k = n - l$.  The base case, $k = 0$, follows from Lemma 
\ref{standardsimple}.  This case also contains the case $n = 2$, so we may
assume that $n > 2$ and $k > 0$.  Proposition \ref{fonlands} now shows
that $$
[\sS_{2n}(\pm(l-1)):L_{2n}(-l)]
=[F' \sS_{2n}(\pm(l-1)):F' L_{2n}(-l)]
= [\sS_{2n-2}(\pm(l-1)):L_{2n-2}(-l)]
,$$ which completes the inductive step.  
The same line of argument proves the
third assertion, namely that $[\sS_{2n}(l-2):L_{2n}(-l)] = 0$ if 
$n \geq 3$ and $l \geq 3$.

Next, we prove that $[\sS_{2n}(\pm(l-1)):L_{2n}(l)] = 0$ by induction on
$k = n - l$.  The case $k = 0$ follows from Lemma \ref{standardsimple}.
If $k > 0$, then Proposition \ref{fonlands} shows
that $$
[\sS_{2n}(\pm(l-1)):L_{2n}(l)]
=[F \sS_{2n}(\pm(l-1)):F L_{2n}(l)]
= [\sS_{2n-2}(\mp(l-1)):L_{2n-2}(-l)]
,$$ which reduces to a previous case and completes the proof of the first
two assertions.
The same line of argument also reduces the
fourth assertion (that $[\sS_{2n}(-l+2):L_{2n}(l)] = 0$ if 
$n \geq 3$ and $l \geq 3$) to previously proved assertions.


%
\end{proof}

The ultimate aim would be to find some ``alcove like'' combinatorics for
the labelling poset for the simple modules. In other words, as in the
Temperley--Lieb  case where the representation theory is controlled by 
$\tilde{\mathrm{A}}_1$ type
alcove combinatorics, we should be able to find a labelling 
poset that gives us alcove combinatorics for some affine Weyl group, 
mirroring the
fact that this algebra is a quotient of a Hecke algebra of
$\tilde{\mathrm{C}}$ type.
Since the above result on composition factors is true for any
specialisation of the parameters, the intriguing consequence of this
result is that the two labels $\pm l$ and $\pm (l-1)$ need never be
comparable (providing $l\ge2$). 

Another consequence of this
proposition using Lemma \ref{lem:noext}
is that there are no non-split extensions between the simple modules
$L_{2n}(\pm l)$ and $L_{2n}(\pm (l-1))$, provided $l \ge 2$.

\begin{lem}
We have $[\sS_{2n}(l-4):L_{2n}(-l)] =0$ for $5 \le l \le n$
and
$[\sS_{2n}(-(l-4):L_{2n}(l)] =0$ for $5 \le l \le n-1$.
\end{lem}
\begin{proof}
If $n > l$, we have $$
[\sS_{2n}(l-4):L_{2n}(-l)]
=[F' \sS_{2n}(l-4):F' L_{2n}(-l)]
= [\sS_{2n-2}(l-4):L_{2n-2}(-l)]
$$ and $$
[\sS_{2n}(-l+4):L_{2n}(l)]
=[F \sS_{2n}(-l+4):F L_{2n}(l)]
= [\sS_{2n-2}(l-4):L_{2n-2}(-l)]
.$$  The result will then follow by induction on $n - l$ if we can show that 
$[\sS_{2n}(n-4):L_{2n}(-n)]$ is zero.


Since all possible composition factors of $\sS_{2n}(n-4)$ (apart from
the simple head $L_{2n}(n-4)$) cannot extend each other, it follows that if
$L_{2n}(-n)$ is a composition factor of $\sS_{2n}(n-4)$ then it must
embed in $\sS_{2n}(n-4)$. Thus it is enough to show that there is no
embedding of $L_{2n}(-n)$ into $\sS_{2n}(n-4)$ when the parameters are
invertible.

Now, $\sS_{2n}(n-4)$ is generated by $ee_2f$ 
(modulo $I_{2n}(n-5)+I_{2n}(-n+5)$) 
and so has basis given by
$$\{ee_1ee_2f, e_1ee_2f, ee_2f, e_3ee_2f,  e_4e_3ee_2f, \ldots,
 e_{n-1}\cdots e_4e_3ee_2f,  fe_{n-1}\cdots  e_4e_3ee_2f \}.$$
Let $v =(a_0, a_1, \ldots, a_n) \in \sS_{2n}(n-4)$ with respect to
 this basis. 
If $v$ generates a one-dimensional
submodule of $\sS_{2n}(n-4)$ isomorphic to the trivial module
($L_{2n}(-n)$) then $e$, $f$, and $e_i$ must all act trivially on
$v$.
Thus
$$0=e v =  (\dl a_0+ a_1,0,\dl a_2, \ldots, \dl a_n) $$
and so $a_2=a_3=\cdots=a_n=0$ as $\dl \ne 0$ 
and $\dl a_0+ a_1=0$.
So $v=( a_0, -\dl a_0, 0,0,\ldots,0)$.
We also need
$$0=f v =  (\dr a_0, -\dr\dl a_0,0,0, \ldots,0)
 $$
and so $a_0=a_1=0$ as $\dr \ne 0$.
Thus $v=0$ and there is no submodule of $\sS_{2n}(n-4)$ isomorphic
to the trivial module.
\end{proof}

We may now produce a poset that works for all parametrisations
for which the parameters are units. 
The next section will show that this poset cannot be coarsened further
for all unit parametrisations.
Thus, we can begin to form a picture of what our
``alcove geometry'' must   look like.

\begin{prop}\label{prop:poset}
The affine symmetric Temperley--Lieb algebra is quasi-hereditary with 
the same standards, $\sS_{2n}(l)$, and  with poset:
$$
\xymatrix@R=5pt{%
& 0 \ar@{-}[dl] \ar@{-}[dd]  \ar@{-}[ddr]  \ar@{-}[drr] \\
-1 \ar@{-}[dd] \ar@{-}[dddr] \ar@{-}[dddrr]& &  & 1 \ar@{-}[dd] \ar@{-}[dddll]
\ar@{-}[dddl]     \\
& -2 \ar@{-}[dd]  \ar@{-}[dddl] \ar@{-}[dddrr]
 & 2 \ar@{-}[dd] \ar@{-}[dddll] \ar@{-}[dddr] &\\
-3 \ar@{-}[dd] \ar@{-}[dddr] \ar@{-}[dddrr]
& & & 3 \ar@{-}[dd] \ar@{-}[dddll]\ar@{-}[dddl] \\
&-4 \ar@{-}[dd]  \ar@{-}[dddl] \ar@{-}[dddrr] 
& 4 \ar@{-}[dd] \ar@{-}[dddll] \ar@{-}[dddr] & \\
-5 \ar@{-}[dd] \ar@{-}[dddr] \ar@{-}[dddrr]  
& & & 5 \ar@{-}[dd] \ar@{-}[dddll]\ar@{-}[dddl]\\
&-6 \ar@{-}[dd]  \ar@{-}[ddl] \ar@{-}[ddrr] 
& 6 \ar@{-}[dd] \ar@{-}[ddll] \ar@{-}[ddr] & \\
-7 \ar@{-}[d] & &
& 7 \ar@{-}[d] \\
{\vdots}&{\vdots}&{\vdots}&{\vdots}
}
$$
\end{prop}

It is possible to draw a more planar version of the above at the cost of
not having
elements that are lower down in the order, lower down on the page:
$$
\epsfbox{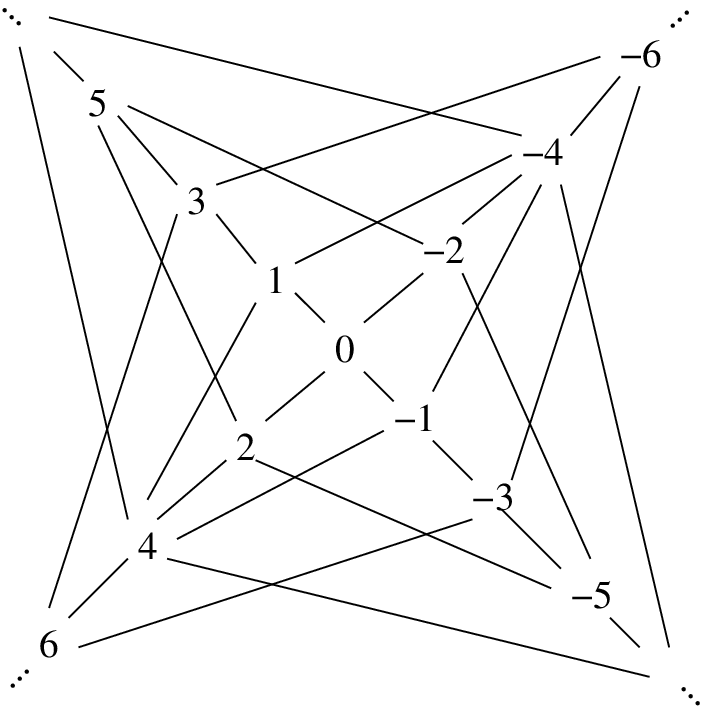}
$$

\section{Some General Gram determinants}\label{sect:gramdets}

The Gram determinant of $\sS_{2n}(0)$ has been determined in
\cite{degiernichols} using a
particular basis of the standard module $\sS_{2n}(0)$, also known as
the ``spine module'' for a particular
parametrisation.
%
In this section we calculate the Gram determinant for the $n+1$
dimensional standard modules,
and give some empirical justification for our chosen parametrisation.
A rigorous definition of a Gram determinant may be found in 
\cite[equation~(40), section~8]{gensymp}. In essence what we do
is define a suitable inner product on the module and then take the determinant
of the square matrix formed by the inner products between all the
basis elements for the module. 

Using the definition of $\sS_{2n}(l)$ we can always find a monomial
basis for $\sS_{2n}(l)$, by letting the generators of $b_{2n}^\phi$
act on $d$, the diagram that generates $\sS_{2n}(l)$.  
If $b_1$, $b_2$ are two such monomial basis elements
then they have the same lower half diagram, which is the same as the
lower half diagram of $d$, denoted $\langle d|$. 
We let $b_2^\circ$ be the
element obtained from $b_2$ by turning all its component diagrams
upside down. When $b_1$ is multiplied by $b_2^\circ$ we obtain a
scalar multiple of $|d\rangle \langle d|$,  possibly zero if $d_2^\circ
d_1$ has fewer than $|2l|$ propagating lines. 
This scalar is then the value of the inner product 
of $b_2$ by $b_1$.  Since for this section $l \ne 0$, we need not
worry about powers of $\kk$ and the Gram determinant may be
calculated in this way. We continue to use $\Deltag_{2n}(l)$
to denote the Gram determinant of $\sS_{2n}(l)$.

\subsection{Some low rank calculations}

We first do some low rank calculations with the original
parametrisation
introduced in section \ref{review}.

\noindent
Case $m=1$:
\\
We see from \cite[figure 3]{gensymp} that 
$b_{2}^{\phi}$ has dimension $5=2^2 +1$, and two standard modules. 
The first one is simple: $\sS_2(-1) = L_2(-1)$ and the second $
\sS_2(0)$,
is two dimensional. Either $\sS_2(0) =L_2(0)$ and
$b_{2}^{\phi}$ is semi-simple or $\dim L_2(0) =1$ and $b_{2}^\phi$
is not semi-simple.
Now $\sS_2(0) \cong b_2^\phi e \cong b_2^\phi  f$,
and so it has basis $\{ e, fe \}$ and Gram matrix
$$
\begin{pmatrix}
\dl & \kk \\
\kk & \dr\kk
\end{pmatrix}
$$
which gives 
the Gram determinant (up to a $\pm$ sign)
\begin{align*}
\Deltag_2(0) &= \kk(\dl \dr - \kk)  \\
&=-\left[{\frac{w_1-w_2+\theta }{2}}\right]
      \left[{\frac{w_1-w_2-\theta }{2}}\right]
\left([w_1][w_2]+\left[{\frac{w_1-w_2+\theta }{2}}\right]
      \left[{\frac{w_1-w_2-\theta }{2}}\right]\right)
\\
&=-\left[{\frac{w_1-w_2+\theta }{2}}\right]
      \left[{\frac{w_1-w_2-\theta }{2}}\right]
([w_1+w_2-1]+[w_1+w_2-3]+\cdots +[w_1-w_2+1]\\
&\hspace{6.5cm}+
[w_1-w_2-1]+[w_1-w_2-3]+\cdots+[\theta+1])
\\
&= - \left[\frac{w_1-w_2+\theta}{2}\right]
\left[\frac{w_1-w_2-\theta}{2}\right]
\left[\frac{w_1+w_2+\theta}{2}\right]
\left[\frac{w_1+w_2-\theta}{2}\right]
\end{align*}
where the last line gives this determinant in terms of the 
parametrisation introduced at the end of subsection \ref{ourpara}.
(NB: to obtain this, we used the quantum number identity
$[a][b]=[a+b-1] +[a+b-3] +\cdots [a-b+1]$.)

\bigskip
\noindent
Case $m=2$:
\\
We see from \cite[figure 3]{gensymp} that 
$\dim b_4^\phi = 10 = 4^2 +1 + 1+1$.
We have $\sS_4(-2)=L_4(-2)$,
$\sS_4(-1)=L_4(-1)$ and $\sS_4(1)=L_4(1)$ as noted in 
Proposition~\ref{standardsimple}.
The dimension of the remaining standard module $\sS_4(0)$ is $4$.
It has Gram determinant (up to a $\pm$ sign)
\begin{align*}
\Deltag_4(0) &=
\kappa_{LR}(\kappa_{LR} - \kappa_L \delta_R )
(\kappa_{LR} - \delta_L \kappa_R )
(\kappa_{LR} - \delta_L \kappa_R - \kappa_L \delta_R
+\delta\kappa_L\kappa_R )\\
&=
\left[\frac{w_1+w_2+\theta+1}{2}\right]
\left[\frac{w_1+w_2-\theta+1}{2}\right]
\left[\frac{w_1-w_2+\theta-1}{2}\right]
\left[\frac{w_1-w_2-\theta-1}{2}\right]\\
& \qquad \times \left[\frac{-w_1+w_2+\theta-1}{2}\right]
\left[\frac{-w_1+w_2-\theta-1}{2}\right]
\left[\frac{w_1+w_2+\theta+3}{2}\right]
\left[\frac{w_1+w_2-\theta+3}{2}\right]
\end{align*}
Here we have suppressed the details of the calculation --- but we note
that the factorising of the determinant was aided by using
globalisation to get the factors of 
$\kappa_{LR} - \kappa_L \delta_R $ and 
$\kappa_{LR} - \delta_L \kappa_R $ and using Gap4 \cite{gap4}.

\subsection{Some more general Gram determinants}
Consider $\sS_{2n}(-(n-2))$ with $n\ge3$.
One basis for $\sS_{2n}(-(n-2))$, generated by $e_1$,
(modulo $I_{2n}(n-3)+I_{2n}(-n+3)$) 
is
$$\{ ee_1,\ e_1,\ e_2e_1,\ e_3e_2e_1,\ e_4e_3e_2e_1,\ \ldots,\
e_{n-1}e_{n-2}\cdots e_1,\ fe_{n-1}e_{n-2}\cdots e_1\}.$$ 
Using this basis the Gram matrix for $\sS_{2n}(-(n-2))$ is
$$ 
\begin{pmatrix}
\dl\kl & \kl &0  &0  &0 &0 &\cdots &0\\
\kl    & \dd &1  &0  &0 &0 &\cdots &0\\
0      & 1   &\dd&1  &0 &0 &\cdots &0\\
0      & 0   & 1 &\dd&1 &0 &\cdots &0\\
\vdots &\vdots &\ddots & \ddots & \ddots & \ddots & \ddots & \vdots \\
0      & 0   &\cdots &0& 1 &\dd&1 &0 \\
0      & 0   &\cdots &0 &0& 1 &\dd&\kr  \\
0      & 0   &\cdots &0 &0& 0 &\kr&\dr\kr  \\
\end{pmatrix}.
$$

Consider $\sS_{2n}(n-3)$ with $n\ge4$.
One basis for $\sS_{2n}(n-3)$, generated by $ee_2$,
(modulo $I_{2n}(n-4)+I_{2n}(-n+4)$) 
is
$$\{ee_1ee_2, e_1ee_2,\ ee_2,\ e_3ee_2,\ e_4e_3ee_2,\ e_4e_3e_2e_1,\ \ldots,\
e_{n-1}e_{n-2}\cdots e_3 e e_2,\ fe_{n-1}e_{n-2}\cdots e_3ee_2\}.$$ 
Using this basis the Gram matrix for $\sS_{2n}(n-3)$ is
$$ 
\begin{pmatrix}
\dl^2\kl&\dl\kl & \dl^2 &0      &0 &0 &\cdots &0\\
\dl\kl  &\dl\dd & \dl   &0      &0 &0 &\cdots &0\\
\dl^2     &\dl    &\dl\dd &\dl  &0 &0 &\cdots &0\\
0       &0      & \dl   & \dl\dd&\dl &0 &\cdots &0\\
\vdots  &\vdots &\ddots & \ddots & \ddots & \ddots & \ddots & \vdots \\
0       &0      &\cdots &0& \dl &\dl\dd&\dl &0 \\
0       &0      &\cdots &0 &0& \dl &\dl\dd&\dl\kr  \\
0       &0      &\cdots &0 &0& 0 &\dl\kr&\dl\dr\kr  \\
\end{pmatrix}.
$$

Consider $\sS_{2n}(-(n-3))$ with $n\ge4$.
One basis for $\sS_{2n}(-(n-3))$, generated by $e_1f$
(modulo $I_{2n}(n-4)+I_{2n}(-n+4)$) 
is
$$\{ ee_1f,\ e_1f,\ e_2e_1f,\ e_3e_2e_1f,\ \ldots,\
e_{n-1}e_{n-2}\cdots e_1f,\ fe_{n-1}e_{n-2}\cdots e_1f\}.$$ 
Using this basis the Gram matrix for $\sS_{2n}(-(n-3))$ is
$$ 
\begin{pmatrix}
\dl\dr\kl&\dr\kl & 0 &0      &0 &0 &\cdots &0\\
\dr\kl  &\dr\dd & \dr   &0      &0 &0 &\cdots &0\\
0     &\dr    &\dr\dd &\dr  &0 &0 &\cdots &0\\
0       &0      & \dr   & \dr\dd&\dr &0 &\cdots &0\\
\vdots  &\vdots &\ddots & \ddots & \ddots & \ddots & \ddots & \vdots \\
0       &0      &\cdots &0& \dr &\dr\dd&\dr &\dr^2 \\
0       &0      &\cdots &0 &0& \dr &\dr\dd&\dr\kr  \\
0       &0      &\cdots &0 &0& \dr^2 &\dr\kr&\dr^2\kr  \\
\end{pmatrix}.
$$

Consider $\sS_{2n}(n-4)$ with $n\ge5$.
We have already seen a basis for this module, generated by $ee_2f$,
in section \ref{fourthree}.
Using this basis the Gram matrix for $\sS_{2n}(n-4)$ is
$$ 
\begin{pmatrix}
\dl^2\dr\kl&\dl\dr\kl & \dl^2\dr &0      &0 &0 &\cdots &0\\
\dl\dr\kl  &\dl\dr\dd & \dl\dr   &0      &0 &0 &\cdots &0\\
\dl^2\dr     &\dl\dr    &\dl\dr\dd &\dl\dr  &0 &0 &\cdots &0\\
0       &0      & \dl\dr   & \dl\dr\dd&\dl\dr &0 &\cdots &0\\
\vdots  &\vdots &\ddots & \ddots & \ddots & \ddots & \ddots & \vdots \\
0       &0      &\cdots &0& \dl\dr &\dl\dr\dd&\dl\dr &\dl\dr^2 \\
0       &0      &\cdots &0 &0& \dl\dr &\dl\dr\dd&\dl\dr\kr  \\
0       &0      &\cdots &0 &0& \dl\dr^2 &\dl\dr\kr&\dl\dr^2\kr  \\
\end{pmatrix}.
$$


\begin{prop}
For $n \ge 3$, we have (up to a $\pm$ sign)
$$
\Deltag_{2n}(-(n-2))= [w_1+1][w_2+1][w_1+w_2-n+2].
$$
Similarly for $n \ge 4$, we have
$$
\Deltag_{2n}(-(n-3))= [w_1]^{n+1}[w_1-1][w_2+1][-w_1+w_2+n-2]
$$ and
$$
\Deltag_{2n}(n-3)= [w_2]^{n+1}[w_2-1][w_1+1][w_1-w_2+n-2],
$$
and for $n\ge 5$, we have
$$
\Deltag_{2n}(n-4)= [w_1]^{n+1}[w_2]^{n+1}[w_1-1][w_2-1][w_1+w_2+n-2].
$$
\end{prop}
\begin{proof}
We prove the fourth statement; the remaining ones are similar.
Removing the factor of $\dl\dr$ from each line and 
expanding the determinant out along the first row we have:
\begin{align*}
\dl^{-n+1}&\dr^{-n+1}\Deltag_{2n}(n-4)\\
&=
\tiny{\left|\begin{matrix}
\dl\kl&\kl & \dl &0      &0 &0 &\cdots &0\\
\kl  &\dd & 1   &0      &0 &0 &\cdots &0\\
\dl  & 1    &\dd & 1 &0 &0 &\cdots &0\\
0       &0      & 1   & \dd& 1&0 &\cdots &0\\
\vdots  &\vdots &\ddots & \ddots & \ddots & \ddots & \ddots & \vdots \\
0       &0      &\cdots &0& 1 &\dd& 1 &\dr \\
0       &0      &\cdots &0 &0& 1 &\dd&\kr  \\
0       &0      &\cdots &0 &0& \dr &\kr&\dr\kr  \\
\end{matrix}\right|}_{n+1} \\
=&\dl\kl 
\tiny{\left|\begin{matrix}
 \dd &1  &0  &0 &0 &\cdots &0\\
 1   &\dd&1  &0 &0 &\cdots &0\\
 0   & 1 &\dd&1 &0 &\cdots &0\\
\vdots &\ddots & \ddots & \ddots & \ddots & \ddots & \vdots \\
 0   &\cdots &0& 1 &\dd&1 &\dr \\
 0   &\cdots &0 &0& 1 &\dd&\kr  \\
 0   &\cdots &0 &0& \dr &\kr&\dr\kr  \\
\end{matrix}\right|}_n
 - \kl^2
\tiny{\left|\begin{matrix}
\dd&1  &0 &0 &\cdots &0\\
 1 &\dd&1 &0 &\cdots &0\\
\vdots & \ddots & \ddots & \ddots & \ddots & \vdots \\
\cdots &0& 1 &\dd&1 &\dr \\
\cdots &0 &0& 1 &\dd&\kr  \\
\cdots &0 &0& \dr &\kr&\dr\kr  \\
\end{matrix}\right|}_{n-1}
\\
 &\qquad\qquad\qquad\qquad\qquad + (2\dl\kl-\dl^2\dd)
\tiny{\left|\begin{matrix}
\dd&1  &0 &0 &\cdots &0\\
 1 &\dd&1 &0 &\cdots &0\\
\vdots & \ddots & \ddots & \ddots & \ddots & \vdots \\
\cdots &0& 1 &\dd&1 &\dr \\
\cdots &0 &0& 1 &\dd&\kr  \\
\cdots &0 &0& \dr &\kr&\dr\kr  \\
\end{matrix}\right|}_{n-2}
\end{align*}
where the subscripts on the square matrices keep note of their size.
Now the determinant of each submatrix is given by
\begin{align*}
&\tiny{\left|\begin{matrix}
 \dd &1  &0  &0 &0 &\cdots &0\\
 1   &\dd&1  &0 &0 &\cdots &0\\
 0   & 1 &\dd&1 &0 &\cdots &0\\
\vdots &\ddots & \ddots & \ddots & \ddots & \ddots & \vdots \\
 0   &\cdots &0& 1 &\dd&1 &\dr \\
 0   &\cdots &0 &0& 1 &\dd&\kr  \\
 0   &\cdots &0 &0& \dr &\kr&\dr\kr  \\
\end{matrix}\right|}_n\\
&=\dr\kr 
\tiny{\left|\begin{matrix}
\dd    &1      &\cdots &0\\
1      &\dd    &\cdots &0\\
\vdots &\ddots & \ddots & \vdots \\
0      &\cdots & \dd & 1  \\
0      &\cdots & 1& \dd   \\
\end{matrix}\right|}_{n-1}
-\kr^2
\tiny{\left|\begin{matrix}
\dd    &1      &\cdots &0\\
1      &\dd    &\cdots &0\\
\vdots &\ddots & \ddots & \vdots \\
0      &\cdots & \dd & 1  \\
0      &\cdots & 1& \dd   \\
\end{matrix}\right|_{n-2}}
+(2 \dr\kr-\dd\dr^2)
\tiny{\left|\begin{matrix}
\dd    &1      &\cdots &0\\
1      &\dd    &\cdots &0\\
\vdots &\ddots & \ddots & \vdots \\
0      &\cdots & \dd & 1  \\
0      &\cdots & 1& \dd   \\
\end{matrix}\right|_{n-3}}\\
&=\dr\kr[n]
-\kr^2[n-1]
 + (2\dr\kr-\dd\dr^2)[n-2]\\
&=[w_2][w_2+1][n]
-[w_2+1]^2[n-1]
 + \Bigl(2[w_2][w_2+1]-[2][w_2]^2\Bigr)[n-2]. \\
\end{align*}
Using the quantum number identities
$$ 
[a+1][b+1] -  [a][b]=[a+b+1]  \qquad \mbox{and} \qquad 
[a+2][b+2] -  [a][b]=[2][a+b+2] $$
we obtain the following for the above determinant:
\begin{align*}
&=[w_2][w_2+1][n]
-[w_2+1]^2[n-1]
 + [w_2]\Bigl([w_2+1]-[w_2-1]\Bigr)[n-2]\\
&=[w_2]\Bigl([w_2+1][n]-[w_2-1][n-2]\Bigr)
-[w_2+1]\Bigl([w_2+1][n-1]-[w_2][n-2]\Bigr)\\
&=[w_2][2][w_2+n-1]
-[w_2+1][w_2+n-1]\\
&=[w_2+n-1][w_2-1].
\end{align*}

Thus
\begin{align*}
[w_1]^{-n+1}[w_2]^{-n+1}&\Deltag_{2n}(n-4)\\
=& \,[w_1][w_1+1][w_2+n-1][w_2-1]
 -[w_1+1]^2[w_2+n-2][w_2-1]\\
&+\Bigl(2[w_1][w_1+1]-[w_1]^2[2]\Bigr)[w_2+n-3][w_2-1]\\
=& \,[w_2-1]\Biggl([w_1][w_1+1][w_2+n-1]
 -[w_1+1]^2[w_2+n-2]\\
&\hspace{4cm}+[w_1]\Bigl([w_1+1]-[w_1-1]\Bigr)[w_2+n-3]\Biggr).\\
\end{align*}
Noting that the expression in large brackets is the same as equation
for the submatrix 
only with $w_2$ replaced with $w_1$ and $n$ replaced
with $w_2+n-1$, we obtain
$$
\Deltag_{2n}(n-4)
=[w_1]^{n-1}[w_2]^{n-1}[w_1-1][w_2-1]
[w_1+w_2+n-2],
$$
as required.
\end{proof}


In order to globalise the factors in these determinants, we need to
know what composition factors correspond to the factors of the Gram
determinant vanishing.
Since $-n+2$, $n-3$, $-n+3$, and $n-4$ are almost minimal in the poset
order,
the only possible other composition factors of these standard modules 
are $L_{2n}(-n)$, $L_{2n}(-n+1)$, $L_{2n}(n-1)$, and $L_{2n}(n-2)$,
the one dimensional modules mentioned in the previous section.
Thus it is possible to find ``by hand'' embeddings of these one
dimensional
modules into the standard ones.

\begin{prop}
We have maps
$$\sS_{2n}(-n) \hookrightarrow \sS_{2n}(-(n-3))\quad \mbox{and} \quad 
\sS_{2n}(-(n-1)) \hookrightarrow \sS_{2n}(n-4)\quad \mbox{for } [w_1-1]=0,$$
$$\sS_{2n}(-n) \hookrightarrow \sS_{2n}(n-3)\quad \mbox{and} \quad
\sS_{2n}(n-1) \hookrightarrow \sS_{2n}(n-4)\quad \mbox{for } [w_2-1]=0,$$
$$\sS_{2n}(-n) \hookrightarrow \sS_{2n}(-(n-2))\quad \mbox{for }
[w_1+w_2-n+2]=0,$$
$$\sS_{2n}(-n+1) \hookrightarrow \sS_{2n}(-(n-3))\quad \mbox{for }
[-w_1+w_2+n-2]=0,$$
$$\sS_{2n}(n-1) \hookrightarrow \sS_{2n}(n-3)\quad \mbox{for }
[w_1-w_2+n-2]=0,$$
$$\sS_{2n}(n-2) \hookrightarrow \sS_{2n}(n-4)\quad \mbox{for }
[w_1+w_2+n-2]=0.$$
\end{prop}
\begin{proof}
In each case, it is a matter of either finding explicit elements of
the larger module that generate a one-dimensional submodule with the
appropriate action (see section \ref{subsect:onedim}), or showing that
the required one-dimensional submodule must exist. 

Consider the module $\sS_{2n}(n-3)$.
We use the basis already introduced for this module.
To show that $\sS_{2n}(-n)$ can be embedded into $\sS_{2n}(n-3)$, we
need a non-zero vector $v\in\sS_{2n}(n-3)$ such that all elements of
$b_{2n}^\phi$ act trivially on it.
Let $v =(a_0, a_1, \ldots, a_n) \in \sS_{2n}(n-3)$ with respect to
the basis given above; note that the coordinates of $v$ are numbered starting
at zero.
We have 
$$e v =  (\dl a_0+ a_1,0,\dl a_2, \ldots, \dl a_n) $$
$$f v =  (0, 0, \ldots,0,a_{n-1}+\dr a_n)$$
$$e_1 v =  (0, \kl a_0+ \dd a_1+a_2,0,0, \ldots, 0) $$
$$e_2 v =  (0,0,\dl a_0+ a_1+\dd a_2+a_3,0,\ldots, 0) $$
$$\mbox{for $3 \le i\le n-2$} \ \ 
  e_i v =  (0,0,\ldots, 0, a_{i-1}+ \dd a_i + a_{i+1},0,\ldots, 0) 
\ \ \mbox{($i$-th position)} $$
$$e_{n-1} v =  (0,0,\ldots, 0, a_{n-2}+ \dd a_{n-1}+\kl a_n,0).$$
It follows that if $ev=0$ then $a_2=a_3=\cdots=a_n=0$, as $\dl \ne 0$ 
and $\dl a_0+ a_1=0$. 
So $v=( a_0, -\dl a_0, 0,0,\ldots,0)$, and thus $f v$ and
$e_i v$ are also zero for $i \ge 2$.
We also need $e_1 v =0$, which gives us that $\kl a_0 +\dd a_1 =0$.
Thus to have a consistent set of equations, we need
$\kl a_0 -\dd\dl a_0 = 0$. So either $a_0=0$ (and hence $v=0$) or 
$\kl -\dd\dl=0=[w_1+1]-[2][w_1]=-[w_1-1]$.
Thus either $[w_1-1] \ne 0$ and $\sS_{2n}(-n)$ does not embed in
  $\sS_{2n}(n-3)$ or $[w_1-1]=0$ and we do have an embedding of
  $\sS_{2n}(-n)$ into $\sS_{2n}(n-3)$.

Similarly, $\sS_{2n}(n-1)$, the module where $b_{2n}^\phi$ acts
trivially except for $e$, which acts as multiplication by $\dl$,
embeds in $\sS_{2n}(n-3)$ if the following set of equations are
satisfied:
\begin{align*}
\dl a_0 + a_1 &=\dl a_0, & a_{n-1} + \dr a_n &=0, \\
 \kl a_0 + \dd a_1 + a_2&=0,
  &\dl a_0 + a_1 + \dd a_2+a_3&=0, \\
   a_{i-1} + \dd a_i+a_{i+1}&=0, \ \mbox{for $3 \le i \le n-2$}, 
& a_{n-2} + \dd a_{n-1} + \kr a_n&=0.
\end{align*}
The system of linear equations in the $a_i$ has $n+1$ by $n+1$ size
coefficient matrix:
$$ \begin{pmatrix}
 0   &1   &0   &0  &0   &0 &0 &\cdots \\
 \kl &\dd &1   &0  &0   &0 &0 &\cdots \\
 \dl &1   &\dd &1  &0   &0 &0 &\cdots \\
 0   &0   &1   &\dd&1   &0 &0 &\cdots \\
 0   &0   &0   &1  &\dd &1 &0 &\cdots \\
\vdots &\ddots & \ddots &\ddots & \ddots & \ddots &\ddots&  \vdots \\
 0   & 0 &\cdots&0 & 1 &\dd&1 &0 \\
 0   & 0 &\cdots&0 &0& 1 &\dd&\kr  \\
 0   & 0 &\cdots&0 &0 &0& 1 &\dr \\
\end{pmatrix}$$
Thus to get non-trivial solutions to the above set of equations
we need the 
matrix to have zero determinant.
This matrix has determinant 
$$-\dr\kl[n-1]
+\dl\dr[n-2]
 + \kl\kr[n-2]
 - \dl\kr[n-3]
=[w_1-w_1+n-2].
$$
Thus $\sS_{2n}(n-1)$ can only embed in $\sS_{2n}(n-3)$ if
$[w_1-w_1+n-2]=0$.

The other maps are constructed in a similar fashion.
\end{proof}

\begin{prop}
We have maps for $n \ge 3$ and $n-m$ even
$$\sS_{2n}(-m) \to \sS_{2n}(-(m-3))\quad \mbox{and} \quad 
\sS_{2n}(-(m-1)) \to \sS_{2n}(m-4)\quad \mbox{for } [w_1-1]=0,$$
$$\sS_{2n}(-m) \to \sS_{2n}(m-3)\quad \mbox{and} \quad
\sS_{2n}(m-1) \to \sS_{2n}(m-4)\quad \mbox{for } [w_2-1]=0,$$
$$\sS_{2n}(-m) \to \sS_{2n}(-(m-2))\quad \mbox{for }
[w_1+w_2-m+2]=0,$$
$$\sS_{2n}(-m+1) \to \sS_{2n}(-(m-3))\quad \mbox{for }
[-w_1+w_2+m-2]=0,$$
$$\sS_{2n}(m-1) \to \sS_{2n}(m-3)\quad \mbox{for }
[w_1-w_2+m-2]=0,$$
$$\sS_{2n}(m-2) \to \sS_{2n}(m-4)\quad \mbox{for }
[w_1+w_2+m-2]=0.$$

We have maps for $n \ge 3$ and $n-m$ odd
$$\sS_{2n}(-m) \to \sS_{2n}(-(m-3))\quad \mbox{for }  [w_1-1]=0,$$
$$\sS_{2n}(-(m-1)) \to \sS_{2n}(m-4)\quad \mbox{for } [w_1-1]=0,$$
$$\sS_{2n}(+m) \to \sS_{2n}(-m+3)\quad \mbox{for } [w_2-1]=0,$$
$$\sS_{2n}(-m+1) \to \sS_{2n}(-m+4)\quad \mbox{for } [w_2-1]=0.$$
\end{prop}
\begin{proof}
To prove this, it is enough to globalise the maps from the previous
proposition, taking care with the swapping of the parameters.
If we assume that $n-m$ is even then we may apply either $G\circ G$ or
$G'\circ G'$ to the maps --- the double application of $G$ or $G'$
ensures
that the swapped parameters are swapped back.

For example: take the map
$\sS_{2n}(-m) \hookrightarrow \sS_{2n}(-(m-3))$
that exists when  $[w_1-1]=0$.
Applying $G'$ once gives a map
$\sS_{2n+2}(-m) \to \sS_{2n+2}(-(m-3))$
for $[w_1+1]=0$. (Note that $G'$ has no effect on $\dl$ and $\kl$ and
hence no effect on $w_1$.)
Applying $G'$ again gives:
$\sS_{2n+4}(-m) \to \sS_{2n+4}(-(m-3))$
for $[w_1+1]=0$. 
Alternatively, we can apply $G\circ G$ to the original
map and get the map above.
\end{proof}

\begin{cor}
The poset in proposition \ref{prop:poset} cannot be coarsened further
and still be a poset for which $b_{2n}^\phi$ is quasi-hereditary for
all parametrisations.
\end{cor}
\begin{proof}
We need to exhibit parametrisations for which each link in the poset
in proposition \ref{prop:poset} is necessary.
To do this it is enough to use the maps from the previous
proposition for links not involving $0$.
For the links involving $0$ we use 
\cite[section 9.2]{gensymp}.
\end{proof}

\section{Quotients of the symplectic blob}


The appearance of $TL_n$ as a subalgebra of $b_n$, 
although constructively natural,
tells us relatively little about their representation theory.
More significant is the appearance of $TL_n$ as a quotient of $b_n$ for
certain special values of the blob parameter \cite{martsaleur}. 
Indeed this is part
of the original paradigm for the generalised alcove geometric approach
\cite{marryom},
which we aim eventually to generalise further to include the symplectic
case. Further, one knows by elementary combinatorial arguments
that $TL_n$ cannot appear as a quotient of $b_n$ ``generically'',
so that, when it does so, we are guaranteed to be studying the
non-generic sector.
By analogy, the study of quotient algebras of the symplectic blob
algebra (with known representation theory)
provides another tool with which to investigate its
representation theory. Again one knows that $TL_n$ and $b_n$ are not generic
quotients, so any such map would embed these known structures into the
non-generic sector.  The study of quotients of the symplectic blob algebra
is more complicated than for the usual blob algebra because of the 
``topological relation''.


\subsection{Some quotients}
\begin{prop}
If $\delta=\kappa_L=\kappa_R$, $\delta_L=\delta_R=\kappa_{LR}=1$
 and $n\ge3$ is odd 
then $b_{n}^x/I \cong TL_{n}$, where $I$ is the ideal generated by
$e-1$ and $f - 1$ and $TL_{n}$ has $q+q^{-1}= \delta$.
\end{prop}
\begin{proof}
Define a map from $b_n^x \to TL_n$ that takes the diagram elements
and removes any blobs. 
Note that this map does not change the underlying structure of any
diagram and lifts to a homomorphism on multiplication of the underlying
diagrams (ignoring the decorations).
This map clearly has kernel equal to $I$ above.

By considering the relations on the diagrams in $b_n^x$ we see that to
make this map a homomorphism from $b_n^x$, 
we need $\delta=\kappa_L=\kappa_R$ and $\delta_L=\delta_R=1$.
(This is by considering the three loop relations and the two blobs = one
blob relations.)
This leaves the ``three blobs'' relation --- which gives us that
$\kappa_{LR}= 1$. As $n$ is odd we do not have the ``topological
relation'' and we do not have the loop with a black and a white blob
replaced by the $\kappa_{LR}$ relation.

The quotient clearly has basis equal to the diagram basis for $TL_n$
and also the same multiplication so we are done.
\end{proof}

\begin{prop}
If $\delta=\kappa_L=\delta_R$, $\delta_L=\kappa_R=\kappa_{LR}=1$
and $n$ is even (and not zero) 
then there is a quotient of $b_{n}^x$ that is isomorphic to
$TL_{n+1}$, where  $q+q^{-1} = \dd$.
\end{prop}
\begin{proof}
We will distinguish the ``$f$'' for $b_{n+1}^x$ with 
a subscript $n+1$.
Consider the composition of maps
$$ b_{n}^x \stackrel{\phi}\to f'b_{n+1}^x f' 
\stackrel{\iota} \hookrightarrow b_{n+1}^x
\stackrel{\psi}\to TL_{n+1},$$
where $\phi$ is the isomorphism of 
(the blob version of) Proposition \ref{fiso}, $f'=\dr^{-1}f_{n+1}$),
$\phi$ is the inverse map to the map $\rho$ of \cite{gensymp} (with no
scalar factor as $\kr =1$), 
$\iota$ is the natural embedding 
and $\psi$ is the quotient map from the previous proposition. 

We claim that $\xi:=\psi\circ \iota \circ \phi$ is onto.
For this, note that $\xi(e_i)$ for $1 \le i \le n-1$ is
the diagram $U_i$ for $TL_{n+1}$.
Also, $\xi(f) = \psi(f_{n+1} e_n f_{n+1}) = U_n$. 
Thus since $\im \xi$ contains all the generators for $TL_{n+1}$, and
it is a homomorphism, it must be onto.
%
%
%
%
\end{proof}

Note that for these propositions, we start with a non-semi-simple 
symplectic blob algebra.  Consider $n$ odd.
Now the $b_{2n}^\phi$ considered has $\dl\dr
-\kk= 1 - 1 = 0$ on one hand. On the other hand, $\dl\dr - \kk =
\left[\frac{w_1+w_2+\theta}{2}\right]
\left[\frac{w_1+w_2-\theta}{2}\right]$. Thus, we have some singular
Gram determinants and at least one of the
standard modules is not simple  and so this
algebra is not semi-simple. 
Hence the quotient, which is isomorphic to the Temperley--Lieb 
algebra, for which we know all decomposition numbers,
gives us information about the possible form of decomposition numbers
of the standard modules for the symplectic blob.

\providecommand{\bysame}{\leavevmode\hbox to3em{\hrulefill}\thinspace}
\providecommand{\MR}{\relax\ifhmode\unskip\space\fi MR }
\providecommand{\MRhref}[2]{%
  \href{http://www.ams.org/mathscinet-getitem?mr=#1}{#2}
}
\providecommand{\href}[2]{#2}

\end{document}